\documentclass[12pt]{article}
\usepackage{latexsym,amsfonts,amssymb}
\usepackage{amsmath}
\usepackage[all]{xy}
\usepackage{hyperref}
\usepackage[dvips]{color}
\usepackage{colordvi,multicol}
\usepackage[marginal]{footmisc}

\def \cal{\mathcal}

\newtheorem{thm}{Theorem}[section]

\newtheorem{lem}[thm]{Lemma}
\newtheorem{pro}[thm]{Proposition}
\newtheorem{defi}[thm]{Definition}

\newtheorem{exa}[thm]{Example}

\numberwithin{equation}{section}
\date{}

\begin{document}
	
	\title{\bf Random periodic solutions of nonautonomous stochastic feedback systems with multiplicative noise}
	\author{}
	\maketitle
	\centerline{Zhao Dong$^{1,2}$, Weili Zhang$^{1,2,*}$ and Zuohuan Zheng $^{3,1,2}$} 
	\centerline{\small $^1$  Academy of Mathematics and Systems Science, Chinese Academy of Sciences, Beijing, 100190, China} 
	\centerline{\small $^2$ School of Mathematical Sciences, University of Chinese Academy of Sciences, Beijing, 100049, China} 
	\centerline{\small $^3$ College of Mathematics and Statistics, Hainan Normal University, Haikou, Hainan 571158, China}
	\centerline{\small dzhao@amt.ac.cn, zhangweili@amss.ac.cn, zhzheng@amt.ac.cn}
	\footnote{* Corresponding author}
	\vskip 1cm  
	\vskip 0.5cm \noindent{\bf Abstract:}\quad 
	We investigate the dynamical behavior of pull-back trajectories for nonautonomous stochastic feedback systems with multiplicative noise. We proved that there exists a random periodic solution of this system and all pull-back trajectories converge to this random periodic solution as time goes to infinitely almost surely. Our results can be applied to nonautonomous stochastic Goodwin negative feedback system, nonautonomous stochastic Othmer-Tyson positive feedback system and nonautonomous stochastic competitive systems etc.
	
	\smallskip
	
	\vskip 0.5cm
	
	
	\vskip 0.3cm
	
	
	\section{Introduction}	
	\quad For deterministic systems, there is a well-developed and constructive theory of systems interconnections and feedback, such as the very successful and widely applied backstepping approach \cite{KKK, I} and stability analysis \cite{DV, G, GJ, HTW, H, S, TO}. Feedback loops play important roles in many biochemical control systems, which often occur in the study of the reaction process in cellular signalling, such as \cite{CD, MCJ, NCB}. It is natural to attempt to extend such work to stochastic systems considering the real word phenomena. As a matter of fact, much excellent research has been done pursuing such extensions, notably studies on stochastic stability \cite{DMS, LJZ, LZJ, M, YXD}. Recently, Marcondes de Freitas and Sontag have initiated a different approach based upon random dynamical systems to investigated the stability of feedback systems involving $real$ $noise$ perturbation in \cite{FS}. Motivated by them, Jiang and Lv considered the global stability of nonlinear stochastic feedback systems driven by additive and multiplicative white noise respectively in \cite{JL, JJL}.
  
   It is natural to attempt to extend nonlinear output(feedback) function to time-preiodic feedback function, that is, consider nonautonomous stochastic systems. Indeed, we have considered the stable random periodic solution of nonautonomous stochastic feedback systems driven by additive white noise in \cite{DZZ}. Our goal in this paper is to prove that there exists a random periodic solutions of nonautonomous stochastic feedback systems driven by multiplicative white noise and all pull-back trajectories converge to this random periodic solution as time tends to infinitely almost surely. In this paper, we will make full use of the theory of random dynamical systems established by L. Arnold \cite{A} and stochastic flows established by H. Kunita \cite{K, KH}. And the powerful theory of monotone random dynamical systems \cite{C} can be applied to investigate the global stability of stochastic flows while the stochastic system admits the stochastic comparison principle, i.e., the system is cooperative or monotone. There are several literatures which use random dynamical systems to investigate the existence of random periodic solution, such as \cite{ZZ, FWZ, FFZ, FZZ}. Compaired with the existing ones, we also consider the global stability.
     
     $$ \xymatrix{& {X_1} \ar[dr]^{}& \\ {X_3} \ar[ur]^{}& & {X_2} \ar[ll]}$$ 
     \vspace{-2em}
     \begin{figure}[ht]
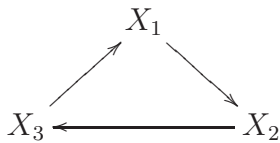
 \caption{Biochemical circuit. The symbol $``X\rightarrow Y"$ means that species $X$ represses the production of species $Y$.} 
     \end{figure}
 
    As a motivation, we first look at a simple biochemical circuit. This biochemical circuit contains three chemical species $X_1$, $X_2$, $X_3$ that interact with one another as shown in Figure 1. Systems of this type are routinely studied as molecular biology, biochemical reaction systems. Furthermore, the strength of the interactions between the species, may depend on enviromental factors such as temperature and the concentrations of other biochemical compounds not explicitly modeled. This dependence may be periodicity and randomicity intrinsically. If this is the case, then a more realistic mathematical model would be a nonautonomous stochastic feedback system of the following form:
    \begin{equation}\label{2221}
    	dx_i=(\alpha_ix_i+h_i(t,x_{i-1}))dt+\sigma_ix_idW_t^i, \ \ i=1,2,3,
    \end{equation}
    which indices taken modulo three, so $x_0=x_3$. Here, $\alpha_i$, $i=1,2,3$ are negative constants, $h_i(\cdot, x_{i-1})$, $i=1,2,3$ are nonincreasing functions in $x_{i-1}$, $h_i(t+T, x_{i-1})=h_i(t, x_{i-1})$, $i=1,2,3$, $T$ is a positive constant, and $W_t(\omega)=\left(W_t^1(\omega),W_t^2(\omega),W_t^3(\omega)\right)$ is a three dimensional standard Brownian motion with $W_0^i(\omega)=0$, $i=1,2,3$, $\omega\in\Omega$.

    \[\xrightarrow{\;v_1\;}X_1\xrightarrow{h_2(\cdot,x_1)}\qquad \xrightarrow{\;v_2\;}X_2\xrightarrow{h_3(\cdot,x_2)}\qquad \xrightarrow{\;v_3\;}X_3\xrightarrow{h_1(\cdot,x_3)}\] 
    \vspace{-2em} 
    \begin{figure}[ht] \caption{Decomposition of the biochemical circuit from Figure 1 into input-output modules. In each partition, $v_i$ indicates the input into the element $X_i$ and $h_i(\cdot,x_{i-1})$ indicates the subsequent output---feedback of the current state.}
    \end{figure}
    
    The nonlinearity of $h_i$ makes the system difficult to study directly.  To overcome this difficulty, we use the decomposition motivated by the work of \cite{FS, JL, JJL}. This idea is to look at (\ref{2221}) as a network of smaller input-output modules as shown in Figure 2, and then we can derive the closed system's properties from these smaller modules' emerging properties. The first step is to open up the feedback loop, rewriting the model as a stochastic system with inputs
    \begin{equation}\label{2222}
    	dx_i=(\alpha_ix_i+v_i(t))dt+\sigma_ix_idW_t^i, \ \ i=1,2,3,
    \end{equation}
    together with a set of outputs 
    \begin{equation}\label{2223}
    	y_i(t)=v_i(t)=h_i(t,x_{i-1}(t)), \ \ i=1,2,3.
    \end{equation}

    Observe that (\ref{2222}) is much easier to study. In fact, in this particular example, we can show that (\ref{2222}) has a unique, globally attracting random periodic solution ${\cal K}(v)$ for each random periodic input $v$. We call ${\cal K}$ defined in this way the input-to-state characteristic of the system. When the open-loop system satisfies certain conditions, the next step is to look at the gain of the system. The output function is read at ${\cal K}(v)$ for each random periodic input $v$, and an operator ${\cal K}^h$ is so defined on the space of random periodic inputs. If ${\cal K}^h$ has a unique, globally attracting random periodic solution, then the input-output system (\ref{2222}) and (\ref{2223}) is said to satisfy the small-gain condition. It is natural to believe that the closed-loop system should have random periodic solution under such circumstances. Periodicity and monotonicity assumptions ensures the system has a unique, globally attracting random periodic solution.
    
    The goal of this paper may be described as to give a rigorous treatment of this example and its generalizations to more general situations. Considering the following $T$-periodic stochastic feedback system with multiplicative linear noise in $\mathbb{R}_+^d$:
    \begin{equation}\label{11}
    	dX_t=\left(AX_t+h(t,X_t)\right)dt+\sum_{k=1}^{d}\sigma_kX_t dW_t^k, 
    \end{equation}       	
    where $W_t=(W_t^1,...,W_t^d)$ is a two-side time Wiener process with values in $\mathbb{R}^d$ on the canonical Wiener space $(\Omega,{\cal F},({\cal F}^t)_{t\in\mathbb{R}},\mathbb{P})$, i.e., ${\cal F}$ is the Borel $\sigma$-algebra of $\Omega=C_0(\mathbb{R},\mathbb{R}^d)=\{\omega:\omega(t)$ continuous, $\omega(0)=\boldsymbol{0},t\in\mathbb{R}\}$; ${\cal F}_s^t$ is the least complete $\sigma$-field for which all $W_u-W_v$, $s\leq v\leq u\leq t$ are measurable and ${\cal F}^t={\cal F}_{-\infty}^t=\bigvee_{s\leq t}{\cal F}_s^t$; $\mathbb{P}$ is the Wiener measure. $A=(a_{ij})_{d\times d}$ is a $(d\times d)$-dimensional matrix. $h:\mathbb{R}\times \mathbb{R}_+^d\rightarrow \mathbb{R}_+^d$, $h(t+T,x)=h(t,x)$ for any $t\in\mathbb{R}, x\in\mathbb{R}_+^d$, $T>0$ is a constant. $\sigma_k$, $k=1, \cdots, d$ are $(d\times d)$-dimensional matrices. 
    
   This paper is organized as follows. In section 2, we review some preliminary concepts and definitions, present the assumptions for the stochastic differential equation (\ref{11}), and define the input-to-state characteristic operator of the system via the pull-back of the discretised stochastic differential equation. In section 3, we describes the asymptotic behavior of stochastic solution, give some auxiliary lemmas, present the definition of gain operator and its properties. In section 4, the main theorem is proved and the global convergence to a unique random periodic solution is presented. In section 5, we present some examples. 

  Convention: Throughout this paper, without loss of generality we always denote a universal set of full $\mathbb{P}$-measure by $\Omega$.

\section{Preliminaries} 
\quad For the convenience of readers, we recall some definitions and basic facts about random dynamical systems and stochastic flows, see \cite{A, C, K, KH} for more details. Let $(\Omega,{\cal F},\mathbb{P})$ be a probability space, $X$ be a Polish space and ${\cal B}(X)$ be its Borel $\sigma$-algebra. Denote $\triangle:=\{(t,s)\in\mathbb{R}^2,s\leq t\}$. 
\begin{defi}[{\cite{A}}]
	A family of mapping on the sample space $\Omega$, $\theta_t:\Omega\rightarrow\Omega$, $t\in \mathbb{R}$ is called a measurable dynamical system if the following conditions are satisfied 	
	\begin{itemize}
		\item [(i)] Identity property: $\theta_0$ is the identity on $\Omega$;
		\item [(ii)] Flow property: $\theta_{t+s}=\theta_{t}\circ \theta_{s}$, where $\circ$ means composition of mappings;
		\item [(iii)] Measurability: $(\omega, t)\mapsto \theta_{t}\omega$ is measurable.
	\end{itemize}
	It is called a measure-preserving or metric dynamical system, if furthermore
	\begin{itemize}
		\item [(iv)] Measure-preserving property: $\mathbb{P}(\theta_{t}(A))=\mathbb{P}(A)$, for every $A\in {\cal F}$ and $t\in \mathbb{R}$.
	\end{itemize}
	In this case, $\mathbb{P}$ is called an invariant measure with respect to the dynamical system $\theta_t$.
\end{defi}
\begin{defi}[{\cite{A}}]
	A (continuous) random dynamical system (RDS) on the Polish space $X$ over a metric dynamical system $(\Omega,{\cal F},\mathbb{P},(\theta_{t})_{t\in\mathbb{R}})$ with time $\mathbb{R}_+$ is a mapping $$\Phi:{\mathbb{R}_+}\times\Omega\times X \rightarrow X, \  (t,\omega,x)\mapsto\Phi(t,\omega,x)$$ 
	which is $({\cal B}({\mathbb R}_+)\otimes {\cal F}\otimes {\cal B}(X), {\cal B}(X))$-measurable and satisfies the following properties:
	\begin{itemize}
		\item [(i)] Continuity: $\Phi(\cdot, \omega, \cdot): \mathbb{R}_+\times X\rightarrow X$, $(t, x)\mapsto \Phi(t, \omega, x)$ is continuous for all $\omega\in\Omega$.
		\item [(ii)] Cocycle property: The mappings $\Phi(t,\omega):=\Phi(t,\omega,\cdot) : X\rightarrow X$ form a cocycle over $\theta(\cdot)$, i.e. they satisfy for all $\omega\in\Omega$,
		\begin{align*}
			\Phi(0,\omega)& \ \text{is the identity on} \ X, \\
			\Phi(t+s,\omega)&=\Phi(t,\theta_s\omega)\circ\Phi(s,\omega),  \ \ \text{for} \ \text{all} \  s,t\in\mathbb{R}_+.
		\end{align*}	
	\end{itemize}
\end{defi}
\begin{defi}[{\cite{A}}]
	A random variable $R: \Omega\rightarrow\mathbb{R}_+$ is called tempered with respect to the dynamical system $\theta$ if
	\begin{equation*}
		\lim\limits_{t\rightarrow \pm\infty}\frac{\log R(\theta_t\omega)}{|t|}=0.
	\end{equation*}	
\end{defi}
This condition is equivalent to the subexponential growth of $t\mapsto R(\theta_t\omega)$,
	\begin{align*}
		\lim\limits_{t\rightarrow \pm\infty}\left\{e^{-\gamma |t|}R(\theta_t\omega)\right\}=0, \ \ \  \text{for any} \ \ \gamma>0,\  \omega\in\Omega,
	\end{align*}
	which implies that, for any $\gamma>0$, $\omega\in\Omega$ 
	\begin{equation*}
		\sup_{t\in\mathbb{R}}\left\{e^{-\gamma |t|}R(\theta_t\omega)\right\}<\infty.
	\end{equation*}
\begin{defi}[{\cite{KH}}]
	A map $\varphi: \triangle\times\Omega\times\mathbb{R}^d\rightarrow\mathbb{R}^d$, $(t,s,x,\omega)\mapsto\varphi(t,s,x,\omega)$ is called a forward stochastic flow if for almost all $\omega\in\Omega$, it satisfies the following properties:
	\begin{itemize}
		\item [(i)]
		$\varphi(t,s,x,\cdot)$ is continuous with respect to $(s, t, x)$.
		\item [(ii)]$\varphi(u,s,\omega)=\varphi(u,t,\omega)\circ\varphi(t,s,\omega)$ holds for all $s\leq t\leq u$.
		\item [(iii)] $\varphi(s,s,\omega)$ is the identity map on $X$ for all $s$.
	\end{itemize}
\end{defi}
\begin{defi}[{\cite{FZZ}}]
	A random periodic solution of period $T>0$ for the forward stochastic flow $\varphi:\triangle\times\Omega\times X\rightarrow X$ is an ${\cal F}$-measurable map $Y: \mathbb{R}\times\Omega\rightarrow X$ such that for almost all $\omega\in\Omega$,
	\begin{equation}\label{14}
		\varphi(t,s,\omega)Y(s,\omega)=Y(t,\omega), \ Y(r+T,\omega)=Y(r,\theta_T\omega), \ \ \text{for any} \ (t, s)\in\triangle, \  r\in\mathbb{R}.
	\end{equation}	
\end{defi}	

Firstly, we consider the corresponding linear homogeneous It$\hat{o}$ stochastic differential equations:
\begin{equation}\nonumber
	dX_t=AX_tdt+\sum_{k=1}^{d}\sigma_kX_tdW_t^k. 
\end{equation}
Without loss of generality, we assume that $\sigma_k$, $k=1, \cdots, d$ has the following form
\begin{align*}
	\sigma_k=\left ( \begin{matrix}
		\sigma_k^1& &\\
		& \ddots &\\
		& & \sigma_k^d	
	\end{matrix}  \right ), \ \ \ \sigma_k^i\in\mathbb{R}, \ \ k, i=1 ,\cdots, d.
\end{align*}
A more general situation can be reduced to this one by a diagonalizing linear transformation.
 According to Definition 3.3.13 in \cite{KS}, it is equivalent to the following Stratonovich stochastic differential equations
\begin{equation}\label{21}
	dX_t=(A-\frac{1}{2}B)X_tdt+\sum_{k=1}^{d}\sigma_kX_t\circ dW_t^k, 
\end{equation}
where 
\begin{align*}
	B=\left ( \begin{matrix}
		\sum_{k=1}^{d}(\sigma_k^1)^2& &\\
		& \ddots &\\
		& & \sum_{k=1}^{d}(\sigma_k^d)^2	
	\end{matrix}  \right ).
\end{align*}

Next, we introduce assumptions that guarantee that the random dynamical system generated by (\ref{21}) in $\mathbb{R}_+^d$ is order-preserving and the existence and uniqueness of solutions for stochastic differential equations (\ref{11}). In order to make use of the technique for monotone systems, we make the following standing assumption on $A$. 
\begin{itemize}
	\item[$\hypertarget{A}{\boldsymbol{(A)}}$] $A$ is cooperative, i.e., $a_{ij}\geq 0$ for all $i,j\in \{1,...,d \}$ and $i\neq j$.
\end{itemize}
Throughout this paper, we will use the norm $|x|:= \max\{|x_i|:i=1,\cdots,d\}$, $|x|_2:= (\sum\limits_{i=1}^d {|x_i|^2})^{\frac{1}{2}}$, $x\in\mathbb{R}^d$ and $||M||_2:= (\sum\limits_{i,j=1}^d {|M_{ij}|^2})^{\frac{1}{2}}$, where $M$ is a $(d\times d)$-dimensional matrix.

Let $\Phi_j(t)=(\Phi_{1j}(t),...,\Phi_{dj}(t)) ^T $ be the solution of equation (\ref{21}) with initial value $X(0)=e_j,j=1,...,d$. Define the $d\times d$ matrix
$$\Phi(t)=(\Phi_1(t),...,\Phi_d(t))=(\Phi_{ij}(t))_{d\times d}.$$
Then $\Phi(t)$ is the fundamental matrix of equation (\ref{21}). It is useful to note that $\Phi(0)$ is the $d\times d$ identity matrix and 
$$d\Phi(t)=A\Phi(t)dt+\sum_{k=1}^{d}\sigma_k\Phi(t)dW_t^k.$$
By Proposition 6.2.2 in \cite{C}, it is clear that (\ref{21}) generates a order-preserving random dynamical system $(\theta, \Phi)$ in $\mathbb{R}_+^d$ and $\Phi(t,\omega)(\mathbb{R}_+^d\setminus \{0\}) \subset \mathbb{R}_+^d\setminus \{0\}$ for any $t\geq 0$, $\omega\in\Omega$, where $\theta$ is the time shift on $\Omega$, i.e.,
\begin{equation}\label{22}
	\theta_t\omega(\cdot):=\omega(t+\cdot)-\omega(t), \ \ t\in\mathbb{R}.
\end{equation}
$\Phi$ satisfies the cocycle property: $\Phi(t+s,\omega)=\Phi(t,\theta_s\omega)\circ\Phi(s,\omega)$ for all $t, s\in\mathbb{R}_+$, $\omega\in\Omega$, and $\Phi(t,\omega)x\geq_{\mathbb{R}_+^d}\Phi(t,\omega)y$ for all $x,y\in\mathbb{R}_+^d$ such that $x\geq_{\mathbb{R}_+^d} y$, where $x \geq_{\mathbb{R}_+^d} y$ means that $x-y\in\mathbb{R}_+^d$. 

Now we give the assumption of $(\theta,\Phi)$ which will be needed in what follows. 
\begin{itemize}
	\item[$\hypertarget{L}{\boldsymbol{(L)}}$] The top Lyapunov exponent for the linear RDS $(\theta,\Phi)$ is a negative real number, i.e., there exist a constant $\lambda>0$ and a stochastic process $R(t,\omega)>0$ which satisfies that for any $\gamma>0$, $\sup_{t\in\mathbb{R}}\left\{e^{-\gamma |t|}\sup_{s\in\mathbb{R}_+}R(s,\theta_t\omega)\right\}<\infty$ such that for all $t\geq 0$, $\omega\in\Omega$,
	\begin{equation}\label{23}
		\parallel \Phi(t,\omega) \parallel:= \max \{|\Phi_{ij}(t,\omega)| :i,j=1,...,d\}\leq R(t,\omega)e^{-\lambda t}.
	\end{equation}
\end{itemize}
\quad In the remainder of this section, we discuss the questions of existence and uniqueness of solutions for stochastic differential equation $(\ref{11})$, as well  as its pull-back trajectories. Let us start with the following assumptions on $h$. 
\begin{itemize}
	\item[$\hypertarget{H}{\boldsymbol{(H)}}$] $h\in C_b^1(\mathbb{R}\times\mathbb{R}_+^d,\mathbb{R}_+^d \setminus \{0\})$, i.e., the function $h$ and its derivatives are both bounded.	And $h$ is order-preserving in $\mathbb{R}_+^d$, i.e., for any $t\in \mathbb{R}$
	$$h(t,x_1)\leq_{\mathbb{R}_+^d} h(t,x_2) \ \ \text{whenever} \ x_1, x_2\in\mathbb{R}_+^d \ \ \text{such that} \ \ x_1\leq_{\mathbb{R}_+^d} x_2,$$
	or anti-order-preserving in $\mathbb{R}_+^d$, i.e.,
	$$h(t,x_1)\geq_{\mathbb{R}_+^d} h(t,x_2) \ \ \text{whenever} \ x_1, x_2\in\mathbb{R}_+^d \ \ \text{such that} \ \ x_1\leq_{\mathbb{R}_+^d} x_2.$$
\end{itemize}

By $\hyperlink{H}{\boldsymbol{(H)}}$, it is easy to check that (\ref{11}) satisfies the conditions of the global Lipschitz and linear growth in $\mathbb{R}_+^d$, since $h$ and its derivatives are both bounded in $\mathbb{R}_+^d$. Let ${\tilde{h}}$ be an extension from $\mathbb{R}\times\mathbb{R}_+^d$ to $\mathbb{R}\times\mathbb{R}^d$ such that ${\tilde{h}}(t,x)=h(t, |x^1|, \cdots, |x^d|)$, $x=(x^1,\cdots, x^d)\in\mathbb{R}^d$. It is clearly that ${\tilde{h}}$ satisfies the conditions of global Lipschitz and linear growth in $\mathbb{R}^d$, we thus have the existence and uniqueness of global solutions for  
\begin{equation}\nonumber
	dX_t=[AX_t+\tilde{h}(t,X_t)]dt+\sum_{k=1}^{d}\sigma_kX_t dW_t^k
\end{equation} 
in $\mathbb{R}^d$, which is a forward stochastic flow $\varphi(t,s,\omega): \mathbb{R}^d\rightarrow \mathbb{R}^d$ (cf.\cite{KS, KH, M, O}). From the form of $\tilde{h}$, we can deduce that the set $\mathbb{R}_+^d$ is forward invariant under the forward stochastic flow, i.e., $\varphi(t,s,\omega)\mathbb{R}_+^d\subset \mathbb{R}_+^d$ for $(t, s)\in\triangle$, $\omega\in\Omega$ and $\varphi(t,s,\omega)x=X(t,s,x,\omega)$, $x\in\mathbb{R}_+^d$ is a unique solution of (\ref{11}). Furthermore,  it can be written as  
\begin{equation}\label{24}
	\varphi(t,s,\omega)x
	=\Phi(t-s,\theta_s\omega)x +\displaystyle \int^{t}_{s}{\Phi(t-r,\theta_r\omega)h(r,\varphi(r,s,\omega)x)dr}, \ (t, s)\in\triangle, \ x\in\mathbb{R}_+^d.
\end{equation} 
 \begin{pro}
 	$\varphi(t,s,\omega)x$ in (\ref{24}) has the following properties:
 	\begin{itemize}
 		\item [(i)] For all $(t, s)\in\triangle$ and $x\in\mathbb{R}_+^d$, $\varphi(t,s,\cdot)x$ is ${\cal F}_s^t$-measurable.
 		\item [(ii)] For all $\omega\in\Omega$, $\varphi(t,s,\omega)x$ is continuous in $(t,s,x)$ and satisfies $\lim\limits_{t\downarrow s}\varphi(t,s,x,\omega)=x$.
 		\item [(iii)] For all $\omega\in\Omega$,
 		\begin{equation}\nonumber
 			\varphi(t,s,\omega)=\varphi(t,r,\omega)\circ\varphi(r,s,\omega), \ \ \text{for all}\ \  s\leq r\leq t, \ \ s,r,t\in\mathbb{R}.
 		\end{equation}
 		\item [(iv)] For all $(t, s)\in\triangle$, $\omega\in\Omega$,
 		\begin{equation}\label{25}
 			\varphi(t+T, s+T, \omega)=\varphi(t, s, \theta_T\omega).
 		\end{equation}
 	\end{itemize}
 \end{pro}
  
  \noindent {\bf Proof.}\ \ The properties (i), (ii), (iii) are obvious. Next we prove (iv). For $(t, s)\in\triangle$, by the periodicity of $h(\cdot, x)$ we have
	\begin{equation}\nonumber
		\begin{aligned}
			\begin{split}
				& \ \  \ \ \varphi(t+T,s+T,\omega)x\\
				&=\Phi(t-s,\theta_{s+T}\omega)x+ \int^{t+T}_{s+T}{\Phi(t+T-r,\theta_r\omega)h(r,\varphi(r,s+T,\omega)x)dr} \\
				&=\Phi(t-s,\theta_{s+T}\omega)x+\int^{t}_{s}{\Phi(t-r,\theta_{r+T}\omega)h(r+T,\varphi(r+T,s+T,\omega)x)dr}\\
				&=\Phi(t-s,\theta_{s+T}\omega)x+\int^{t}_{s}{\Phi(t-r,\theta_{r+T}\omega)h(r,\varphi(r+T,s+T,\omega)x)dr}
			\end{split}
		\end{aligned}
	\end{equation}
	We thus have found that the function 	  
	$$\psi(r,s,\theta_T\omega)x=\varphi(r+T,s+T,\omega)x, \ \ s\leq r $$
	satisfies
	$$\psi(t,s,\theta_T\omega)x=\Phi(t-s,\theta_{s+T}\omega)x+\int^{t}_{s}{\Phi(t-r,\theta_{r+T}\omega)h(r,\psi(r,s,\theta_T\omega)x)dr}. $$	  
	By the uniqueness of the solution
	$$\varphi(t,s,\theta_T\omega)x=\psi(t,s,\theta_T\omega)x=\varphi(t+T,s+T,\omega)x.$$
	which implies that (\ref{25}) holds.
\hfill\fbox\\ 

 For any $n\in\mathbb{N}_+$, $t\geq -nT$, let $\varphi(t,-nT,\omega)x$ be the solution with the initial value $X(-nT)=x$, by (\ref{24}) we have
\begin{equation}\nonumber
	\varphi(t,-nT,\omega)x
	=\Phi(t+nT,\theta_{-nT}\omega)x + \displaystyle \int^{t}_{-nT}{\Phi(t-s,\theta_s\omega)h(s,\varphi(s,-nT,\omega)x)ds}.
\end{equation}  
By $\hyperlink{L}{\boldsymbol{(L)}}$, it is evident that for all $x\in\mathbb{R}_+^d$, $\omega\in\Omega$, $\lim_{t\rightarrow\infty}\Phi(t, \theta_{-t}\omega)x=0$. Regarding the feedback function $h$ as an input term, we can define the input-to-state characteristic operator $\cal K$ associated with given inputs in $\mathbb{R}_+^d$ as follows:
\begin{equation}\label{26}
	\cal{K}(v)(t,\omega)=\displaystyle \int^{t}_{-\infty}{\Phi(t-s,\theta_s\omega)v(s,\omega)ds}, \ \ t\in\mathbb{R},\ \omega\in\Omega,
\end{equation}
where the stochastic process $v:\mathbb{R}\times\Omega\rightarrow\mathbb{R}_+^d$ is bounded. 
	It is evident that ${\cal K}$ is well defined by condition $\hyperlink{L}{\boldsymbol{(L)}}$. In fact, by (\ref{23}) we have $||\Phi(t,\omega)||_2\leq d||\Phi(t,\omega)||\leq d R(t,\omega)e^{-\lambda t}$, $\lambda>0$, and so for any $t\in\mathbb{R}, \omega\in\Omega$, 
	\begin{align*}
		&~~~~\displaystyle \int^{t}_{-\infty}|{\Phi(t-s,\theta_s\omega)v(s,\omega)|_2ds}\\
		&\leq d \displaystyle\int^{t}_{-\infty}R(t-s,\theta_s\omega)e^{-\lambda (t-s)}|v(s+t,\omega)|_2ds\\
		&\leq Cd 
		\sup_{s\in\mathbb{R}} \big\{e^{-\frac{1}{2}\lambda|s|}\sup_{r\in\mathbb{R}_+}R(r,\theta_s\omega)\big\} \displaystyle\int^{t}_{-\infty}e^{-\lambda (t-s)+\frac{\lambda}{2}|s|}ds\\
		&<\infty.
	\end{align*} 

\section{Measurability and asymptotic behavior}
\quad In this section, we give some lemmas to describe the measurability and the dynamical behavior of the pull-back trajectory which will be used in the proof of our main result. Throughout this paper, we set $$\boldsymbol{N}=(N_1,...,N_d),\ \ N_i= \sup_{t\in\mathbb{R}, x\in\mathbb{R}_+^d}|h_i(t,x)|, \ \ i=1,...,d.$$
	
	\begin{lem}\label{3.1}
		For each $n\in\mathbb{N}^+$, $x\in\mathbb{R}_+^d$, $\alpha\in [0, T)$, let
		$$a_n^h(t, \omega)=\inf{\{h(t, \varphi(t,t-\alpha-mT,\omega)x):m\geq n,m\in\mathbb{N}_+\}},\ \ t\in\mathbb{R}, \omega\in\Omega,$$
		and
		$$b_n^h(t, \omega)=\ \sup{\{h(t, \varphi(t,t-\alpha-mT, \omega)x):m\geq n, m\in\mathbb{N}_+\}},\ \ t\in\mathbb{R}, \omega\in\Omega,$$
		where $\inf$ and $\sup$ mean the greatest lower bound and the least upper bound, respectively. Then $a_n^h(t, \omega)$ and $b_n^h(t, \omega)$ are progressively measurable with respect to $\{{\cal F}^t\}$.	
	\end{lem}
	\noindent {\bf Proof.}\ \ We only prove the case of $b_n^h(t, \omega)$ for the sake of convenience and the case of $a_n^h(t, \omega)$ can be proved analogously.
	 Define
	\begin{equation*}
		B_n^h(t, \omega):={\{h(t, \varphi(t,t-\alpha-mT,\omega)x):m\geq n,m\in\mathbb{N}_+\}}.
	\end{equation*}
	First we show that $b_n^h(t, \omega)$ are well defined. For each $n\in\mathbb{N}^+$, $x\in\mathbb{R}_+^d$, by the boundedness of $h$ we know that $B_n^h(t, \omega)$ is a bounded set for any $t\in\mathbb{R}$, $\omega\in\Omega$, which implies that $B_n^h(t, \omega)$ is order-bounded. And since $\mathbb{R}_+^d$ is strongly minihedral (see \cite{C}, Definition 3.1.7), $b_n^h(t, \omega)$ exists.	Next, we prove that $b_n^h(t, \omega)$ are progressively measurable with respect to $\{{\cal F}^t\}$. Define
	\begin{align*}
		\beta_{n,M}^h(t,\omega)=\sup\Big\{h(t, \varphi(t,t-\alpha-mT,\omega)x):n\leq m\leq M, m\in\mathbb{N}_+\Big\}.
	\end{align*}
	By the continuity of $h$, the properties of $\varphi$ and Corollary 3.1.1(ii) in \cite{C}, $\beta_{n,M}^h(t,\omega)$ is progressively measurable with respect to $\{{\cal F}^t\}$ for every $M=1,2,\cdots$. It is clear that
	\begin{equation*}
		\beta_{n,1}^h(t,\omega)\leq \beta_{n,2}^h(t,\omega)\leq \cdots\leq \beta_{n,M}^h(t,\omega)\leq\cdots .
	\end{equation*}
	Moreover, by the boundedness of $h$ in $\mathbb{R}_+^d$, $b_n^h(t, \omega)=\lim\limits_{{M} \to {\infty}}\beta_{n,M}^h(t,\omega)$ is a progressively measurable with respect to $\{{\cal F}^t\}$.  
	\hfill\fbox\\
	
	\begin{lem}
		Assume that conditions $\hyperlink{A}{\boldsymbol{(A)}}$, $\hyperlink{L}{\boldsymbol{(L)}}$ and $\hyperlink{H}{\boldsymbol{(H)}}$ hold. Then we have
		\begin{equation}\label{31}
			{\cal K}(\underline{\lim}h(\cdot, \varphi))\leq  \underline{\lim}\varphi\leq  \overline{\lim}\varphi\leq {\cal K}( \overline{\lim}h(\cdot, \varphi)),
		\end{equation}  
		where 
		$$\underline{\lim}\varphi(t,\omega):=\lim \limits_{{n} \to {\infty}} \inf\{\varphi(t, t-\alpha-mT, \omega)x:m\geq n, m\in\mathbb{N}_+\}, \ \  t\in\mathbb{R}, \omega\in\Omega,$$
		
		$$\overline{\lim}\varphi(t,\omega):=\lim \limits_{{n} \to {\infty}}\sup\{\varphi(t, t-\alpha-mT, \omega)x:m\geq n, m\in\mathbb{N}_+\}, \ \ t\in\mathbb{R}, \omega\in\Omega,$$
		
		$$\underline{\lim}h(\cdot, \varphi)(t,\omega):=\lim \limits_{{n} \to {\infty}} a_n^h(t,\omega), \ \ t\in\mathbb{R}, \omega\in\Omega,$$

		$$\overline{\lim}h(\cdot, \varphi)(t,\omega):=\lim \limits_{{n} \to {\infty}}b_n^h(t,\omega), \ \ t\in\mathbb{R}, \omega\in\Omega,$$
		for each $x\in\mathbb{R}_+^d$, $\alpha\in [0, T)$.
	\end{lem}
	
	\noindent {\bf Proof.}\ \ By the definition of $\underline{\lim}\varphi$ and $\overline{\lim}\varphi$ we have $\underline{\lim}\varphi\leq  \overline{\lim}\varphi$, the second inequality in (\ref{31}) is proved.  For the sake of convenience we only prove the first inequality in (\ref{31}) and the third inequality can be proved analogously. Similar to Lemma 3.1, we can easily get that $ \underline{\lim}\varphi(t,\omega)$ and $ \underline{\lim}h(\cdot, \varphi)(t,\omega)$ exist, which are also progressively measurable with respect to $\{{\cal F}^t\}$. Then by the boundedness of $h$, (\ref{26}) and the Fubini Theorem, ${\cal K}( \underline{\lim}h(\cdot, \varphi))(t,\omega)$ is well defined, progressively measurable with respect to $\{{\cal F}^t\}$. By (\ref{26}), the definition of the $ \underline{\lim}h(\cdot, \varphi)$ and Lebesgue's dominated convergence theorem, we have ${\cal K}( \underline{\lim}h(\cdot, \varphi))=\lim\limits_{{n} \to {\infty}}{\cal K}(a_n^h)$. For fixed $n\in\mathbb{N}_+$, it is enough to prove that 
	\begin{equation}\nonumber
		\begin{aligned}
			\begin{split}
				&\ \ \ \ {\cal K}(a_n^h)(t, \omega)\\
				& = \displaystyle \int^{t}_{-\infty}{\Phi(t-s, \theta_s\omega)\inf\{h(s,\varphi(s, s-\alpha-mT, \omega)x):m\geq n, m\in\mathbb{N}_+\}ds}\\
				&=\lim \limits_{{\tilde{m}} \to {\infty}\atop{\tilde{m}}\geq n} \Big\{\Phi(t+{\tilde{m}}T, \theta_{-{\tilde{m}}T}\omega)x \\
				&~~~~+\displaystyle \int^{t}_{n T-{\tilde{m}T}}{\Phi(t-s, \theta_s\omega)\inf\{h(s,\varphi(s, s-\alpha-mT, \omega)x):m\geq n, m\in\mathbb{N}_+\}ds}\Big\}\\
				&=\lim \limits_{\tilde{n} \to {\infty}\atop \tilde{n}\geq n}\inf \Big\{\Phi(t+{\tilde{m}}T, \theta_{-{\tilde{m}}T}\omega)x \\
				&+\displaystyle \int^{t}_{nT-{\tilde{m}T}}{\Phi(t-s, \theta_s\omega)\inf\{h(s,\varphi(s, s-\alpha-mT, \omega)x):m\geq n, m\in\mathbb{N}_+\}ds}:{\tilde{m}}\geq \tilde{n}, {\tilde{m}}\in\mathbb{N}_+\Big\}\\
				&\leq \lim \limits_{\tilde{n} \to {\infty}\atop \tilde{n}\geq n}\inf \Big\{\Phi(t+{\tilde{m}}T, \theta_{-{\tilde{m}}T}\omega)x \\
				&~~~~+\displaystyle \int^{t}_{n T-{\tilde{m}T}}{\Phi(t-s, \theta_s\omega)h(s,\varphi(s, s-\alpha-{\tilde{m}}T, \omega)x)ds}:{\tilde{m}}\geq \tilde{n},{\tilde{m}}\in\mathbb{N}_+\Big\}\\    
				&\leq \lim \limits_{\tilde{n} \to {\infty}}\inf \Big\{\Phi(t+{\tilde{m}}T, \theta_{-{\tilde{m}}T}\omega)x \\
				&~~~~+\displaystyle \int^{t}_{-{\tilde{m}T}}{\Phi(t-s, \theta_s\omega)h(s,\varphi(s, s-\alpha-{\tilde{m}}T, \omega)x)ds}:{\tilde{m}}\geq \tilde{n},{\tilde{m}}\in\mathbb{N}_+\Big\}\\ 
				&= \underline{\lim}\varphi(t,\omega),
			\end{split}
		\end{aligned}
	\end{equation}
	where the third equality has used Lemma A.2 in \cite{FS}, while the second-to-last inequality has applied the order-preserving property of $\Phi$ and the positivity of $h$.
	\hfill\fbox\\   
	
	\begin{lem}
		Assume that conditions $\hyperlink{A}{\boldsymbol{(A)}}$, $\hyperlink{L}{\boldsymbol{(L)}}$ and $\hyperlink{H}{\boldsymbol{(H)}}$ hold. Then we have the following:\\
		$(i)$ if $h$ is order-preserving in $\mathbb{R}_+^d$, then 
		\begin{equation}\label{32}
			h(\cdot, \underline{\lim}\varphi)\leq  \underline{\lim}h(\cdot,\varphi)\leq  \overline{\lim}h(\cdot,\varphi)\leq h(\cdot,  \overline{\lim}\varphi);
		\end{equation}  
		$(ii)$ if $h$ is anti-order-preserving in $\mathbb{R}_+^d$, then 
		\begin{equation}\label{33}
			h(\cdot,  \overline{\lim}\varphi)\leq  \underline{\lim}h(\cdot,\varphi)\leq  \overline{\lim}h(\cdot,\varphi)\leq  h(\cdot, \underline{\lim}\varphi).
		\end{equation} 
	\end{lem}
	
	\noindent {\bf Proof.}\ \ Indeed, the proof of the first inequality in (\ref{32}) is adequate and the rest of the results of this lemma can be obtained analogously. Observe that $h$ is order-preserving in $\mathbb{R}_+^d$, then for fixed $n,k\in \mathbb{N}_+, k\geq n$,
	$$h(t,\inf\{\varphi(t, t-\alpha-mT, \omega)x:m\geq n, m\in\mathbb{N}_+\})\leq h(t,\varphi(t, t-\alpha-kT, \omega)x),$$
thus we have
	\begin{align}\label{34}
		&~~~~h(t,\inf\{\varphi(t, t-\alpha-mT, \omega)x:m\geq n,m\in\mathbb{N}_+\})\nonumber\\
		&\leq \inf \{h(t,\varphi(t, t-\alpha-mT, \omega)x):m\geq n,m\in\mathbb{N}_+\}. 
	\end{align} 
 By the definition of $\underline{\lim}\varphi$, the continuity of $h$, (\ref{34}) and the definition of $\underline{\lim}h(\cdot, \varphi)$, we have
	\begin{equation}\nonumber
		\begin{aligned}
			\begin{split}
				h(\cdot, \underline{\lim}\varphi)(t,\omega)
				&= h(t,\lim \limits_{{n} \to {\infty}}\inf\{\varphi(t, t-\alpha-mT, \omega)x:m\geq n,m\in\mathbb{N}_+\})\\
				&= \lim \limits_{{n} \to {\infty}}h(t,\inf\{\varphi(t, t-\alpha-mT, \omega)x:m\geq n,m\in\mathbb{N}_+\})\\
				&\leq  \lim \limits_{{n} \to {\infty}} \inf\{h(t,\varphi(t, t-\alpha-mT, \omega)x):m\geq n,m\in\mathbb{N}_+\}\\
				&=\underline{\lim}h(\cdot, \varphi)(t,\omega).
			\end{split}
		\end{aligned}
	\end{equation}   
The proof is complete.
	\hfill\fbox\\  
	
	\begin{lem}
		Assume that conditions $\hyperlink{A}{\boldsymbol{(A)}}$, $\hyperlink{L}{\boldsymbol{(L)}}$ and $\hyperlink{H}{\boldsymbol{(H)}}$ hold. We have
		\begin{equation}\label{35}
			{\cal K}(a_n^h)\leq \underline{\lim}\varphi\leq \overline{\lim}\varphi\leq{\cal K}(b_n^h), \ \ n\in\mathbb{N}_+,
		\end{equation}    
		where $a_n^h$ and $b_n^h$ are as defined in Lemma 3.1. Furthermore, we define the gain operator
		$${\cal K}^h(u)(t, \omega)=h(t,{\cal K}(u)(t, \omega)),\ \ t\in\mathbb{R},\ \omega\in\Omega,$$ 
		and we have the following:\\
		$(i)$ if $h$ is order-preserving in $\mathbb{R}_+^d$, then for fixed $n\in\mathbb{N}_+$,
		\begin{equation}\label{36}
			({\cal K}^h)^k(a_n^h)\leq  \underline{\lim}h(\cdot,\varphi)\leq  \overline{\lim}h(\cdot,\varphi)\leq ({\cal K}^h)^k(b_n^h), \ \ k\in\mathbb{N}_+.
		\end{equation}  
		$(ii)$ if $h$ is anti-order-preserving in $\mathbb{R}_+^d$, then for fixed $n\in\mathbb{N}_+$,
		\begin{equation}\label{37}
			({\cal K}^h)^{2k}(a_n^h)\leq  \underline{\lim}h(\cdot,\varphi)\leq  \overline{\lim}h(\cdot,\varphi)\leq ({\cal K}^h)^{2k}(b_n^h), \ \ k\in\mathbb{N}_+.
		\end{equation} 
	\end{lem}
	
	\noindent {\bf Proof.}\ \  By the definition of $a_n^h$ and $b_n^h$, it is evident that 
	$$a_n^h\leq  \underline{\lim}h(\cdot,\varphi)\leq  \overline{\lim}h(\cdot,\varphi)\leq b_n^h,\ \ n\in\mathbb{N}_+.$$
	By the order-preserving property of $\Phi$, ${\cal K}(u)$ is monotone with respect to $u$, and consequently   
	$${\cal K}(a_n^h)\leq {\cal K}( \underline{\lim}h(\cdot,\varphi))\leq {\cal K}(  \overline{\lim}h(\cdot,\varphi))\leq {\cal K}(b_n^h), \ \  n\in\mathbb{N}_+.$$
	By (\ref{31}), we have
	$${\cal K}(a_n^h)\leq  \underline{\lim}\varphi\leq  \overline{\lim}\varphi\leq {\cal K}(b_n^h),  \ \  n\in\mathbb{N}_+.$$
	which implies that (\ref{35}) holds. 
	
	In what follows, we claim that (\ref{36}) and (\ref{37}) hold.
	
	If $h$ is order-preserving in $\mathbb{R}_+^d$, then it deduces that $h$ preserves the inequalities in (\ref{35}):
	$${\cal K}^h(a_n^h)\leq h(\cdot, \underline{\lim}\varphi)\leq h(\cdot, \overline{\lim}\varphi)\leq {\cal K}^h(b_n^h),  \ \  n\in\mathbb{N}_+.$$ 
	which together with (\ref{32}) implies
	$${\cal K}^h(a_n^h)\leq  \underline{\lim}h(\cdot,\varphi)\leq  \overline{\lim}h(\cdot,\varphi)\leq {\cal K}^h(b_n^h),  \ \  n\in\mathbb{N}_+.$$  
	This proves (\ref{36}) for $k=1$. Next we assume that, for some $k\in\mathbb{N}$,
	$$({\cal K}^h)^k(a_n^h)\leq  \underline{\lim}h(\cdot,\varphi)\leq  \overline{\lim}h(\cdot,\varphi)\leq ({\cal K}^h)^k(b_n^h),  \ \  n\in\mathbb{N}_+,$$
	holds. From the monotonicity of ${\cal K}$ and (\ref{31}), we have:
	\begin{equation}\nonumber
		\begin{aligned}
			\begin{split}
				{\cal K}\left(({\cal K}^h)^k(a_n^h)\right)
				&\leq {\cal K}( \underline{\lim}h(\cdot,\varphi))\leq  \underline{\lim}\varphi\\
				&\leq  \overline{\lim}\varphi\leq {\cal K}( \overline{\lim}h(\cdot,\varphi))\leq  {\cal K}\left(({\cal K}^h)^k(b_n^h)\right).
			\end{split}
		\end{aligned}
	\end{equation} 
	By the monotonicity of $h$ in $\mathbb{R}_+^d$ and (\ref{32}), we get that
	$$({\cal K}^h)^{k+1}(a_n^h)\leq  \underline{\lim}h(\cdot,\varphi)\leq  \overline{\lim}h(\cdot,\varphi)\leq ({\cal K}^h)^{k+1}(b_n^h),  \ \  n\in\mathbb{N}_+.$$ 
	Therefore, we conclude that (\ref{36}) holds by mathematical induction.
	
	If $h$ is anti-order-preserving in $\mathbb{R}_+^d$, similar to $h$ is order-preserving in $\mathbb{R}_+^d$, we deduce that
	$${\cal K}^h(b_n^h)\leq h(\cdot, \overline{\lim}\varphi)\leq h(\cdot, \underline{\lim}\varphi)\leq {\cal K}^h(a_n^h),  \ \  n\in\mathbb{N}_+.$$
	by (\ref{33}), we have
	$${\cal K}^h(b_n^h)\leq \underline{\lim}h(\cdot,\varphi)\leq  \overline{\lim}h(\cdot,\varphi)\leq {\cal K}^h(a_n^h),  \ \  n\in\mathbb{N}_+.$$
	Combining the monotonicity of ${\cal K}$ and (\ref{31}), it shows that
	$${\cal K}\left({\cal K}^h(b_n^h)\right)\leq \underline{\lim}\varphi\leq  \overline{\lim}\varphi\leq {\cal K}\left({\cal K}^h(a_n^h)\right), \ \  n\in\mathbb{N}_+,$$
	which together with the anti-monotonicity of $h$ in $\mathbb{R}_+^d$ and (\ref{33}) implies
	$$({\cal K}^h)^2(a_n^h)\leq  \underline{\lim}h(\cdot,\varphi)\leq  \overline{\lim}h(\cdot,\varphi)\leq ({\cal K}^h)^2(b_n^h),  \ \  n\in\mathbb{N}_+.$$ 
	The rest of the proof of (\ref{37}) can be obtained analogously to $h$ is order-preserving in $\mathbb{R}_+^d$ by the mathematical induction.
	\hfill\fbox\\

	\section{ Main results} 
	\quad In this section, we state our main result on the stable random periodic solution of nonautonomous stochastic feedback system (\ref{11}) and prove them. 
	
	Let ${\cal M}$ be the space of all progressively measurable (with respect to $\{{\cal F}^t\}$) processes
	$f:\mathbb{R}\times\Omega\rightarrow [\boldsymbol{0},\boldsymbol{N}]$, and 
	$$f(t+T, \omega)=f(t, \theta_T\omega), \ \ \text{for any} \ \ t\in\mathbb{R},\ \omega\in\Omega.$$ 
	
	Here the metric on ${\cal M}$ is given as follows:
	$$\rho(f_1, f_2):=\mathop{ \sup}\limits_{t\in\mathbb{R}}\mathbb{E}{\left|f_1(t,\omega)-f_2(t,\omega)\right|},  \ \ \text{for all}\ f_1, f_2\in{\cal M}.$$
	
	\begin{lem}
		$({\cal M}, \rho)$ is a complete metric space.
	\end{lem}
	\noindent {\bf Proof.}\ \ It is clear that $({\cal M}, \rho)$ is a metric space. To prove completeness assume that $\{f_n, n\in\mathbb{N}\}$ is a Cauchy sequence  in $({\cal M}, \rho)$, i.e.
	\begin{equation*}
		\mathop{ \sup}\limits_{t\in\mathbb{R}}\mathbb{E}{|f_n(t,\omega)-f_m(t,\omega)|}\rightarrow 0, \ \ \text{as}\ \ m,n\rightarrow \infty.
	\end{equation*}
	Then one can find a subsequence $\{f_{n_k}\}$ such that 
	\begin{equation*}
		\mathop{ \sup}\limits_{t\in\mathbb{R}}\mathbb{E}{|f_{n_{k+1}}(t,\omega)-f_{n_{k}}(t,\omega)|}\leq 2^{-k}.
	\end{equation*}
	Thus
	\begin{align*}
		&~~~~\sup_{t\in\mathbb{R}}(\mathbb{E}{\sum_{k=1}^{\infty}|f_{n_{k+1}}(t,\omega)-f_{n_{k}}(t,\omega)|}+\mathbb{E}|f_{n_1}(t,\omega)|)\\
		&\leq\sum_{k=1}^{\infty}\sup_{t\in\mathbb{R}}\mathbb{E}{|f_{n_{k+1}}(t,\omega)-f_{n_{k}}(t,\omega)|}+\sup_{t\in\mathbb{R}}\mathbb{E}|f_{n_1}(t,\omega)|<\infty, 
	\end{align*}
	which implies that for any $t\in\mathbb{R}$,  
	\begin{equation*}
		{\cal N}_t:=\left\{\omega\in \Omega ~\Big|~\sum_{k=1}^{\infty}|f_{n_{k+1}}(t,\omega)-f_{n_{k}}(t,\omega)|+|f_{n_1}(t,\omega)|=+\infty\right\}
	\end{equation*}
	is a null set and thus
	\begin{equation*}
		\sum_{k=1}^{\infty}|f_{n_{k+1}}(t,\omega)-f_{n_{k}}(t,\omega)|+|f_{n_1}(t,\omega)|<\infty,\ \ t\in\mathbb{R}, \ \omega\in {\cal N}_t^c.
	\end{equation*}
	So for any $t\in\mathbb{R}, \omega\in {\cal N}_t^c$, $\lim\limits_{k\rightarrow\infty}f_{n_k}(t, \omega)$ exists and $$\lim\limits_{k\rightarrow\infty}f_{n_k}(t, \omega)= \sum_{k=1}^{\infty}(f_{n_{k+1}}(t,\omega)-f_{n_{k}}(t,\omega))+f_{n_1}(t,\omega).$$ 
	For any $t\in\mathbb{R}, \omega\in\Omega$, set $$f^j(t, \omega)=\liminf_{k\rightarrow\infty}f_{n_k}^j(t, \omega),\ \ j=1,\cdots,d$$ 
	and 
	$$f(t,\omega)=(f^1(t, \omega), \cdots, f^d(t, \omega)).$$
	It is clear that $f$ is well defined, progressively measurable (with respect to $\{{\cal F}^t\}$), and $f\in [\boldsymbol{0},\boldsymbol{N}]$. Since
	\begin{equation*}
		f_{n_k}(t+T, \omega)=f_{n_k}(t, \theta_T\omega), \ \ \text{for any} \ \ t\in\mathbb{R},\ \omega\in\Omega.
	\end{equation*}
	We have 
	\begin{equation*}
		f^j(t+T, \omega)=\liminf_{k\rightarrow\infty}f_{n_k}^j(t+T, \omega)=\liminf_{k\rightarrow\infty} f_{n_k}^j(t, \theta_T\omega)= f^j(t, \theta_T\omega), 
	\end{equation*}
	which implies that
	\begin{equation*}
		f(t+T, \omega)=f(t, \theta_T\omega), \ \ \text{for any} \ \ t\in\mathbb{R},\ \omega\in\Omega.
	\end{equation*}
	Therefore, $f\in{\cal M}$.
	
	In fact, $f(t,\omega)=\sum_{l=k}^{\infty}(f_{n_{l+1}}(t,\omega)-f_{n_{l}}(t,\omega))+f_{n_k}(t,\omega)$ for any $t\in\mathbb{R}$, $\omega\in\ {\cal N}_t^c$. And thus 
	\begin{align*}
		\sup_{t\in\mathbb{R}}
		\mathbb{E}|f_{n_k}-f|
		=\sup_{t\in\mathbb{R}}\mathbb{E}\left|\sum_{l=k}^{\infty}(f_{n_{l+1}}-f_{n_l})\right|
		\leq \sum_{l=k}^{\infty}2^{-l}=2^{1-k},
	\end{align*} 
	which implies that $d(f, f_{n_k})\rightarrow 0$ as $k\rightarrow \infty$. Hence $d(f, f_n)\rightarrow 0$ as $k\rightarrow \infty$ since $\{f_n\}$ is a Cauchy sequence. Thus $({\cal M}, \rho)$ is a complete metric space.
	\hfill\fbox\\
	
	\begin{lem}
		Assume that conditions $\hyperlink{A}{\boldsymbol{(A)}}$, $\hyperlink{L}{\boldsymbol{(L)}}$, $\hyperlink{H}{\boldsymbol{(H)}}$ hold. Assume additionally that the following conditions on $R$ and $h$ are satisfied:
		\begin{itemize}
			\item[$\hypertarget{R}{\boldsymbol{(R)}}$] Let $R(t-s,\theta_s\omega)$, $(t,s)\in\triangle$ be ${\cal F}_s^t$- measurable and $\sup_{t\in\mathbb{R}_+}\mathbb{E}(R(t,\omega))<\infty$.
			\item[$\hypertarget{H_1}{\boldsymbol{(H_1)}}$] Let $L:=  \max \{\ \sup_{t\in\mathbb{R},x\in\mathbb{R}_+^d} |\frac{\partial h_i(t,x)}{\partial x_j}|, i,j=1,...,d\}$ such that $\sup_{t\in\mathbb{R}_+}\mathbb{E}(R(t,\omega))<\frac{\lambda}{Ld^2}$. 
		\end{itemize}
	    
	    Then the gain operator ${\cal K}^h: {\cal M}\rightarrow {\cal M}$, $f\mapsto{\cal K}^h(f) $ is a contractive mapping, where ${\cal K}^h(f)(t, \omega)=h(t,[{\cal K}(f)](t, \omega))$. 
	\end{lem}
	\noindent {\bf Proof.}\ \ First, we prove that  ${\cal K}^h:{\cal M}\rightarrow{\cal M}$ is well defined. For any $f\in ({\cal M}, \rho)$, by the definition of ${\cal K}^h$, ${\cal K}^h(f)\in [\boldsymbol{0}, \boldsymbol{N}]$. And by (\ref{26}), the order-preserving of $\Phi$ and the positivity of $f$, for any $t\in\mathbb{R},\ \omega\in\Omega$, we have
	\begin{equation}\nonumber
		\begin{aligned}
			\begin{split}
				[{\cal K}(f)](t,\theta_T\omega)
				& = \displaystyle \int^{t}_{-\infty}{\Phi(t-r,\theta_{r+T}\omega)f(r,\theta_T\omega)dr}\\
				&=\displaystyle \int^{t+T}_{-\infty}{\Phi(t+T-r,\theta_r\omega)f(r-T,\theta_T\omega)dr}\\
				&= \displaystyle \int^{t+T}_{-\infty}{\Phi(t+T-r,\theta_r\omega)f(r,\omega)dr}\\
				&= [{\cal K}(f)](t+T,\omega).
			\end{split}
		\end{aligned}
	\end{equation}
	By the fact that $h(t+T, x)= h(t, x), t\in\mathbb{R}, x\in\mathbb{R}_+^d$, we have
	\begin{equation}\nonumber
		\begin{aligned}
			\begin{split}
				{\cal K}^h(f)(t+T, \omega)
				& = h(t+T,{\cal K}(f)(t+T, \omega))\\
				&=h(t,{\cal K}(f)(t, \theta_T\omega))\\
				&= {\cal K}^h(f)(t, \theta_T\omega).
			\end{split}
		\end{aligned}
	\end{equation}
	By $\hyperlink{H}{\boldsymbol{(H)}}$, the measurability of $\Phi$, and the Fubini theorem, it is evident that ${\cal K}^h(f)$ is progressively measurable with respect to $\{{\cal F}^t\}$. So ${\cal K}^h:{\cal M}\rightarrow{\cal M}$ is well defined. 
	
	Next we prove that ${\cal K}^h$ is a contractive mapping. By $\hyperlink{H_1}{\boldsymbol{(H_1)}}$, we have
	\begin{equation}\label{441}
	\mathop{ \sup}\limits_{t\in\mathbb{R}, x\in\mathbb{R}_+^d}{||\triangledown_x\left(h(t,x)\right)||}=\mathop{ \sup}\limits_{t\in\mathbb{R}, x\in\mathbb{R}_+^d}{\Big|\Big|\Big(\frac{\partial h_i(t,x)}{\partial x_j}\Big)_{d\times d}\Big|\Big|}\leq L.
	\end{equation} 
	For any $f_1, f_2 \in ({\cal M}, \rho)$, by the fact that $|\Phi(x)|\leq d||\Phi||\cdot|x|$ for all $x\in\mathbb{R}_+^d$, $\Phi\in\mathbb{R}^{d\times d}$, (\ref{441}) and $\hyperlink{L}{\boldsymbol{(L)}}$, we get
	\begin{equation}\nonumber
		\begin{aligned}
			\begin{split}
				&~~~~\rho({\cal K}^h(f_1), {\cal K}^h(f_2))\\
				&=\mathop{ \sup}\limits_{t\in\mathbb{R}}\mathbb{E}\left|{\cal K}^h(f_1)-{\cal K}^h(f_2)\right|\\
				& = \mathop{ \sup}\limits_{t\in\mathbb{R}}\mathbb{E}\left|h(\cdot,{\cal K}(f_1))-h(\cdot,{\cal K}(f_2))\right|\\
				& = \mathop{ \sup}\limits_{t\in\mathbb{R}}\mathbb{E}\left|\int_{0}^{1}{\triangledown_xh(\cdot,{\cal K}(f_2)+r({\cal K}(f_1)-{\cal K}(f_2)))}dr\cdot[{\cal K}(f_1)-{\cal K}(f_2)]\right|\\
				&\leq  d\mathop{ \sup}\limits_{t\in\mathbb{R}, x\in\mathbb{R}_+^d}{||\triangledown_x(h(t,x))||}\cdot\mathop{ \sup}\limits_{t\in\mathbb{R}}\mathbb{E}|{\cal K}(f_1)-{\cal K}(f_2)|\\
				&\leq Ld\mathop{ \sup}\limits_{t\in\mathbb{R}}\mathbb{E}\left|\displaystyle \int^{t}_{-\infty}{\Phi(t-s,\theta_s\omega)f_1(s,\omega)}ds-\displaystyle \int^{t}_{-\infty}{\Phi(t-s,\theta_s\omega)f_2(s,\omega)}ds\right|\\
				&\leq Ld^2\mathop{ \sup}\limits_{t\in\mathbb{R}}\mathbb{E}\displaystyle \int^{t}_{-\infty}{||\Phi(t-s,\theta_s\omega)||\cdot|f_1(s,\omega)-f_2(s,\omega)|}ds\\
				&\leq Ld^2\mathop{ \sup}\limits_{t\in\mathbb{R}}\mathbb{E}\displaystyle \int^{t}_{-\infty}{e^{-\lambda(t-s)}R(t-s,\theta_s\omega)\cdot|f_1(s,\omega)-f_2(s,\omega)|}ds\\
				&\leq Ld^2\mathop{ \sup}\limits_{t\in\mathbb{R}}\displaystyle \int^{t}_{-\infty}e^{-\lambda(t-s)}\displaystyle \int_{\Omega}{R(t-s,\theta_s\omega)\cdot|f_1(s,\omega)-f_2(s,\omega)|}P(d\omega)ds\\
				&= Ld^2\mathop{ \sup}\limits_{t\in\mathbb{R}}\displaystyle \int^{t}_{-\infty}{e^{-\lambda(t-s)}\mathbb{E}(R(t-s,\theta_s\omega))\mathbb{E}|f_1(s,\omega)-f_2(s,\omega)|}ds\\
				&=Ld^2\mathop{ \sup}\limits_{t\in\mathbb{R}}\displaystyle \int^{t}_{-\infty}{e^{-\lambda(t-s)}\mathbb{E}(R(t-s,\omega))\mathbb{E}|f_1(s,\omega)-f_2(s,\omega)|}ds\\
				&\leq Ld^2\mathop{ \sup}\limits_{t\in\mathbb{R}_+}\mathbb{E}(R(t,\omega))\mathop{ \sup}\limits_{t\in\mathbb{R}}\mathbb{E}|f_1(t,\omega)-f_2(t,\omega)|\mathop{ \sup}\limits_{t\in\mathbb{R}}\displaystyle \int^{t}_{-\infty}{e^{-\lambda(t-s)}}ds\\
				&=\frac{Ld^2\mathop{ \sup}\limits_{t\in\mathbb{R}_+}\mathbb{E}(R(t,\omega))}{\lambda}\rho(f_1, f_2)\\
				&<\rho(f_1, f_2),
			\end{split}
		\end{aligned}
	\end{equation}
	where the third-to-last equality holds because of $\hyperlink{R}{\boldsymbol{(R)}}$ and the independence of $R(t-s,\theta_s\omega)$ and $f_1(s,\omega)-f_2(s,\omega)$, the third-to-last inequality holds because of the $\mathbb{P}$-measure preserving property of $\theta$.
	\hfill\fbox\\
	
	\begin{thm}\label{4.3}
		Assume that conditions $\hyperlink{A}{\boldsymbol{(A)}}$, $\hyperlink{L}{\boldsymbol{(L)}}$, $\hyperlink{H}{\boldsymbol{(H)}}$, $\hyperlink{R}{\boldsymbol{(R)}}$ and $\hyperlink{H_1}{\boldsymbol{(H_1)}}$ hold. Then there exist a unique nonnegative fixed point $u\in ({\cal M}, \rho)$ for the gain operator ${\cal K}^h$ such that for any $t\in \mathbb{R}$, 
		\begin{equation*}
		     \lim \limits_{{n} \to {\infty}}{\varphi(t, -nT, \cdot)x}={\cal K}(u)(t, \cdot),\ \  x\in\mathbb{R}_+^d, \  \mathbb{P}-a.s.
		\end{equation*}
	 	Furthermore, for almost all $\omega\in\Omega$,  
		\begin{equation}\nonumber
			\varphi(t,s,\omega){\cal K}(u)(s,\omega)={\cal K}(u)(t,\omega),\  {\cal K}(u)(r+T,\omega)={\cal K}(u)(r,\theta_T\omega), \ \text{for any} \ (t, s)\in\triangle,\ r\in\mathbb{R}.
		\end{equation}	
		i.e., ${\cal K}(u)$ is a random periodic solution of the forward stochastic flow generated by (\ref{11}) in $\mathbb{R}_+^d$ and for any $t\in \mathbb{R}$, all pull-back trajectories originating from nonnegative orthant converge to this positive random periodic solution almost surely. 
	\end{thm}
	\noindent {\bf Proof.}\ \  By (\ref{36}), (\ref{37}), regardless of the monotonicity or anti-monotonicity for $h$, for fixed $n\in\mathbb{N}_+$, we have
	\begin{equation}\label{42}
		({\cal K}^h)^{2k}(a_n^h)\leq \underline\lim h(\cdot, \varphi)\leq\overline\lim h(\cdot, \varphi)\leq ({\cal K}^h)^{2k}(b_n^h), \ \ k\in\mathbb{N}_+
	\end{equation} 
	where $a_n^h$ and $b_n^h$ are as defined in Lemma \ref{3.1}. By Lemma \ref{3.1}, $a_n^h$ and $b_n^h$ are bounded progressively measurable processes with respect to $\{{\cal F}^t\}$. By the definition of $a_n^h$, $h(t+T, x)= h(t, x), t\in\mathbb{R}, x\in\mathbb{R}_+^d$ and (\ref{25}), for any $t\in\mathbb{R}$, $\omega\in\Omega$, we have  
	\begin{equation}\nonumber
		\begin{aligned}
			\begin{split}
				a_n^h(t+T, \omega)
				& = \inf{\{h(t+T, \varphi(t+T, t+T-\alpha-mT,\omega)x):m\geq n,m\in\mathbb{N}_+\}}\\
				& = \inf{\{h(t, \varphi(t, t-\alpha-mT,\theta_T\omega)x):m\geq n,m\in\mathbb{N}_+\}}\\
				&= a_n^h(t, \theta_T\omega).
			\end{split}
		\end{aligned}
	\end{equation}
	Similarly, we have $b_n^h(t+T, \omega)= b_n^h(t, \theta_T\omega)$ for any $t\in\mathbb{R}$, $\omega\in\Omega$. And it is clear that $a_n^h, b_n^h\in [\boldsymbol{0},\boldsymbol{N}]$. Therefore, $a_n^h$ and $b_n^h$ are both in $({\cal M}, \rho)$. 
	
	Since ${\cal K}^h$ is a contractive mapping on the complete metric space $({\cal M}, \rho)$, by the Banach fixed point theorem (\cite{Y}), there exists a unique $u\in {\cal M}$  such that for any $t\in\mathbb{R}$,
	$$[{\cal K}^h(u)](t, \omega)=u(t, \omega), \ \ \mathbb{P}-a.s.$$
	and
	\begin{equation}\nonumber
		\lim \limits_{{k} \to {\infty}}\sup_{t\in\mathbb{R}}\mathbb{E}|[({\cal K}^h)^{2k}(a_n^h)]-u|= \lim \limits_{{k} \to {\infty}}\sup_{t\in\mathbb{R}}\mathbb{E}|[({\cal K}^h)^{2k}(b_n^h)]-u|=0,
	\end{equation}
	which implies that $[({\cal K}^h)^{2k}(a_n^h)](t, \cdot)\overset{\mathbb{P}}{\rightarrow}u(t, \cdot)$ and $[({\cal K}^h)^{2k}(b_n^h)](t, \cdot)\overset{\mathbb{P}}{\rightarrow}u(t, \cdot)$ for all $t\in\mathbb{R}$. Therefore, there exists a subsequence $\{k_j\}_{j\in\mathbb{N}}$ such that for any $t\in\mathbb{R}$,
	\begin{equation}\label{43}
		\lim \limits_{{j} \to {\infty}}{[({\cal K}^h)^{2k_j}(a_n^h)](t, \cdot)}=u(t, \cdot)= \lim \limits_{{j} \to {\infty}}{[({\cal K}^h)^{2k_j}(b_n^h)](t, \cdot)}, \ \ \mathbb{P}-a.s.
	\end{equation}
	Combining (\ref{42}) and (\ref{43}), we obtain that for any $t\in\mathbb{R}$,
	$$[ \underline\lim h(\cdot, \varphi)](t, \cdot)=[ \overline\lim h(\cdot, \varphi)](t, \cdot)=u(t, \cdot) \ \ \mathbb{P}-a.s.$$
	which together with (\ref{31}) implies that for any $t\in\mathbb{R}$,
	\begin{equation}\label{44}
		[ \underline\lim\varphi](t, \cdot)=[ \overline\lim \varphi](t, \cdot)={\cal K}(u)(t, \cdot) \ \ \mathbb{P}-a.s.
	\end{equation}
	By the definition of infimum and supremum, it is clear that
	\begin{equation}\nonumber
		\begin{aligned}
			\begin{split}
				&~~~~\inf\{\varphi(t, t-\alpha-mT, \omega)x:m\geq n, m\in\mathbb{N}_+\}\\
				&\leq {\varphi(t,t-\alpha-nT, \omega)x}\\
				&\leq \ \sup\{\varphi(t, t-\alpha-mT, \omega)x:m\geq n, m\in\mathbb{N}_+\}, \ \ x\in\mathbb{R}_+^d.
			\end{split}
		\end{aligned}
	\end{equation}
	Let $n\rightarrow\infty$ in the above inequality, by (\ref{44}), for any $t\in\mathbb{R}$,  we have
	\begin{equation}\label{45}
		[ \underline\lim\varphi](t, \omega)=[ \overline\lim \varphi](t, \omega)=\lim \limits_{{n} \to {\infty}}\varphi(t, t-\alpha-nT, \omega)x, \ \ \mathbb{P}-a.s.
	\end{equation}
	Combining (\ref{44}) and (\ref{45}), we obtain that for any $t\in\mathbb{R}$, $t=m_0T+\alpha_0, m_0\in\mathbb{Z}, \alpha_0\in[0,T)$,
	$$\lim \limits_{{n} \to {\infty}}{\varphi(t, -nT, \cdot)x}=\lim \limits_{{n} \to {\infty}}{\varphi(t, t-\alpha_0-nT, \cdot)x}={\cal K}(u)(t, \cdot),\ \ x\in\mathbb{R}_+^d,\  \mathbb{P}-a.s.$$
	which proves that for any $t\in \mathbb{R}$, 
	\begin{equation}\label{46}
		\lim \limits_{{n} \to {\infty}}{\varphi(t, -nT, \cdot)x}={\cal K}(u)(t, \cdot), \ \ x\in\mathbb{R}_+^d, \ \mathbb{P}-a.s..
	\end{equation} 
	
	By (\ref{26}), the progressively measurability with respect to $\{{\cal F}^t\}$ of $u$ and the Fubini Theorem, ${\cal K}(u)(t,\omega)$ is continuous and progressively measurable with respect to $\{{\cal F}^t\}$.
	
	By the continuity of $\varphi$ in $\mathbb{R}_+^d$ and (\ref{46}), we can show that for any fixed $(t, s)\in\triangle$,
	\begin{align}\label{47}
		\varphi(t,s,\omega){\cal K}(u)(s,\omega)
		& = \varphi(t,s,\omega) \lim \limits_{{n} \to {\infty}}{\varphi(s,-nT, \omega)x}\nonumber\\
		& = \lim \limits_{{n} \to {\infty}}\varphi(t,s,\omega){\varphi(s,-nT, \omega)x}\nonumber\\
		&=\lim \limits_{{n} \to {\infty}}{\varphi(t, -nT, \omega)x}\nonumber\\
		&={\cal K}(u)(t, \omega),\ \ \mathbb{P}-a.s.
	\end{align}
	
	Next, we claim that for almost all $\omega\in\Omega$, $$\varphi(t,s,\omega){\cal K}(u)(s,\omega)={\cal K}(u)(t, \omega), \ \ (t, s)\in\triangle.$$ By the continuity of $\varphi(t,s,\omega,x)$ in $(t,s,x)$ for all $\omega\in\Omega$ and the continuity of ${\cal K}(u)(t,\omega)$ in $t$ for all $\omega\in\Omega$, we have
	\begin{align*}
		{\cal N}_{s,t}&:=\left\{\omega\in\Omega:\varphi(t,s,\omega){\cal K}(u)(s,\omega)\not={\cal K}(u)(t, \omega)\right\}\\
		&=\bigcup_{k=1}^{\infty} \left\{\omega\in\Omega: |\varphi(t,s,\omega){\cal K}(u)(s,\omega)-{\cal K}(u)(t, \omega)|>\frac{1}{k} \right\}\\
		&\subset \bigcup_{k=1}^{\infty} \left\{\omega\in\Omega: |\varphi(q_t,p_s,\omega){\cal K}(u)(p_s,\omega)-{\cal K}(u)(q_t, \omega)|>\frac{1}{2k} \right\}\\
		&\subset \bigcup_{k=1}^{\infty}\bigcup_{(q,p)\in \mathbb{Q}^2\cap \triangle}\left\{\omega\in\Omega: |\varphi(q,p,\omega){\cal K}(u)(p,\omega)-{\cal K}(u)(q, \omega)|>\frac{1}{2k} \right\}\\
		&=\bigcup_{(q,p)\in \mathbb{Q}^2\cap \triangle} \bigcup_{k=1}^{\infty} \left \{\omega\in\Omega: |\varphi(q,p,\omega){\cal K}(u)(p,\omega)-{\cal K}(u)(q, \omega)|>\frac{1}{2k} \right\}\\
		&=\bigcup_{(q,p)\in \mathbb{Q}^2\cap \triangle}{\cal N}_{p,q}.
	\end{align*}
	Thus we obtain $\mathbb{P}(\bigcup_{(t,s)\in \triangle}{\cal N}_{s,t})\leq \mathbb{P}(\bigcup_{(q,p)\in \mathbb{Q}^2\cap \triangle}{\cal N}_{p,q})=0$, since $\mathbb{P}({\cal N}_{s,t})=0, (t,s)\in\triangle$ (by (\ref{47})).
	
	Furthermore, $u(t+T, \omega)=u(t, \theta_T\omega)$, for any $t\in\mathbb{R}, \omega\in\Omega$ since $u\in {\cal M}$. Thus by (\ref{26}) we have
	\begin{align*}
		{\cal K}(u)(t,\theta_T\omega)
		& = \displaystyle \int^{t}_{-\infty}{\Phi(t-r,\theta_{r+T}\omega)u(r,\theta_T\omega)dr}\\
		&=\displaystyle \int^{t+T}_{-\infty}{\Phi(t+T-r,\theta_r\omega)u(r-T,\theta_T\omega)dr}\\
		&= \displaystyle \int^{t+T}_{-\infty}{\Phi(t+T-r,\theta_r\omega)u(r,\omega)dr}\\
		&= {\cal K}(u)(t+T,\omega), \ \ \ t\in\mathbb{R},\ \ \omega\in\Omega.
	\end{align*}
	So for almost all $\omega\in\Omega$,
		\begin{equation*}
		\varphi(t,s,\omega){\cal K}(u)(s,\omega)={\cal K}(u)(t,\omega),\  {\cal K}(u)(r+T,\omega)={\cal K}(u)(r,\theta_T\omega), \ \ \text{for any} \ (t, s)\in\triangle,\ r\in\mathbb{R}.
	\end{equation*}	
	\hfill\fbox\\

	\section{Examples}
	\quad In this section, we show the efficiency of our result. Our main results Theorem \ref{4.3} works for the  $T$-periodic stochastic Goodwin negative feedback system, $T$-periodic stochastic Othmer-Tyson positive feedback system and $T$-periodic stochastic competitive systems. Our main task is to check the condition $\hyperlink{L}{\boldsymbol{(L)}}$, $\hyperlink{R}{\boldsymbol{(R)}}$ and $\hyperlink{H_1}{\boldsymbol{(H_1)}}$ in order to use Theorem \ref{4.3}. That is to say, we need to choose a suitable $\lambda>0$ and random process $R$.
	
	Now we consider nonautonomous stochastic single loop feedback system
	\begin{eqnarray}\label{51}
		\left\{
		\begin{array}{l}
			dx_1=\left(-\alpha_1 x_1+ f(t,x_n)\right)dt+\sigma_1x_1dW_t^1,\\
			dx_i=(x_{i-1}-\alpha_i x_i)dt+\sigma_ix_idW_t^i, \ 2\leq i\leq n,
		\end{array}
		\right.
	\end{eqnarray}
	where $\alpha_i>0$, $\sigma_i>0$ for $i=1,\cdots,n$ and $f\in C_b^1(\mathbb{R}\times\mathbb{R}_+, \mathbb{R}_+)$, i.e., $f$ and its derivatives are both bounded. Moreover, we assume that $f(t+T, x)=f(t, x)$, $t\in\mathbb{R}$, $x\in\mathbb{R_+}$ and $f(t, \cdot)$, $t\in\mathbb{R}$ is increasing or decreasing in $\mathbb{R}_+$. 
	
	The corresponding linear homogeneous stochastic It$\hat{o}$ type differential equations is
	\begin{eqnarray}\label{52}
		\left\{
		\begin{array}{l}
			dx_1=-\alpha_1 x_1dt+\sigma_1x_1dW_t^1,\\
			dx_i=(x_{i-1}-\alpha_i x_i)dt+\sigma_ix_idW_t^i, \ 2\leq i\leq n.
		\end{array}
		\right.
	\end{eqnarray}
	By the variation-of-constants formula, we can easily calculate the fundamental matrix $\Phi(t,\omega)$ of (\ref{52}) as follows:
	$$
	\Phi(t, \omega)=\begin{pmatrix}
		\Phi_{11}(t,\omega) & 0 & \cdots & 0 \\
		\Phi_{21}(t,\omega) & \Phi_{22}(t,\omega) & \cdots & 0 \\
		\vdots & \vdots & \ddots & \vdots \\
		\Phi_{n1}(t,\omega) & \Phi_{n2}(t,\omega) & \cdots & \Phi_{nn}(t,\omega) \\ 
	\end{pmatrix}
	$$
	for all $t\geq 0$ and $\omega\in\Omega$, where
	\begin{equation}\label{53}
		\Phi_{ij}(t, \omega)=\left\{
		\begin{array}{rcl}
			\displaystyle \int^{t}_{0}{\Phi_{ii}(t-s,\theta_{s}\omega)\Phi_{i-1, j}(s,\omega)ds}, & & {1\leq j \leq i-1}\\ 
			e^{-(\alpha_i+\frac{1}{2}\sigma_i^2)t+\sigma_i W_t^i(\omega)},& & {j=i}\\
			0, & & {i+1\leq j\leq n}
		\end{array} \right.
	\end{equation}
	for all $i=1, \cdots, n$. Let $\lambda=\frac{1}{n+1}\min\{\alpha_1, \cdots, \alpha_n\}$, then it is easy to check that
	\begin{align}\label{54}
		\Phi_{ii}(t,\omega)&\leq e^{-[(n+1)\lambda+\frac{1}{2}\sigma_i^2]t+\sigma_i W_t^i(\omega)}\nonumber\\
		&=e^{-(i\lambda+\frac{1}{2}\sigma_i^2)t+\sigma_i W_t^i(\omega)}e^{-(n+1-i)\lambda t}\nonumber\\
		&=R_{ii}(t, \omega)e^{-(n+1-i)\lambda t},
	\end{align}
	for all $t\geq 0$ and $\omega\in\Omega$, where
	\begin{equation}\label{55}
		R_{ii}(t, \omega)=e^{-(i\lambda+\frac{1}{2}\sigma_i^2)t+\sigma_i W_t^i(\omega)}, \ \ i=1, \cdots, n.
	\end{equation}
	We prove that $\sup_{s\in\mathbb{R}_+}R_{ii}(s,\omega)$ is a tempered random variable, for any $\gamma>0$, by (\ref{55}) we have the following estimate
	\begin{align*}
		&~~~~\sup_{t\in\mathbb{R}}\big\{e^{-\gamma |t|}\sup_{s\in\mathbb{R_+}}R_{ii}(s,\theta_t\omega)\big\}\\
		&=\sup_{t\in\mathbb{R}}\Big\{e^{-\gamma |t|}\sup_{s\in\mathbb{R_+}}\big\{ e^{-(i\lambda+\frac{1}{2}\sigma_i^2)s+\sigma_i W_s^i(\theta_t\omega)}\big\}\Big\}\\
		&=\sup_{s\in\mathbb{R_+},t\in\mathbb{R}}\Big\{e^{-\gamma|t|-(i\lambda+\frac{1}{2}\sigma_i^2)s+\sigma_i \big(W_{s+t}^i(\omega)-W_t^i(\omega)\big)}\Big\}\\
		&\leq \sup_{s\in\mathbb{R_+},t\in\mathbb{R}}\Big\{e^{-\left(\frac{\gamma}{2}\wedge \left(i\lambda+\frac{1}{2}\sigma_i^2\right)\right)\left|t+s\right|+\sigma_i W_{s+t}^i(\omega)}e^{-\frac{\gamma}{2}|t|-\sigma_i W_t^i(\omega)}\Big\}\\
		&\leq \sup_{t\in\mathbb{R}}\Big\{e^{-\left(\frac{\gamma}{2}\wedge \left(i\lambda+\frac{1}{2}\sigma_i^2\right)\right)\left|t\right|+\sigma_i W_{t}^i(\omega)}\Big\}\sup_{t\in\mathbb{R}}\Big\{e^{-\frac{\gamma}{2}|t|-\sigma_i W_t^i(\omega)}\Big\}\\
		&< \infty,
	\end{align*} 
	where the last inequality holds because of the law of the iterated logarithm of Brownian motions.
	
	Next, we claim that 
	\begin{equation}\label{56}
		\Phi_{ij}(t,\omega)\leq e^{-(n+1-i)\lambda t}R_{ij}(t,\omega), \ \ 1\leq j\leq i-1
	\end{equation} 
	where 
	\begin{equation}\label{57}
		R_{ij}(t,\omega)=\int^{t}_{0}e^{-\lambda s}{R_{ii}(t-s,\theta_{s}\omega)R_{i-1, j}(s,\omega)ds},  \ \ 1\leq j\leq i-1.
	\end{equation}
	In order to check (\ref{56}), we only present the proof of $\Phi_{21}(t,\omega)$ and $\Phi_{31}(t,\omega)$, the rest can be analogously completed by induction. Combining (\ref{53}), (\ref{54}), we have
	\begin{align*}
		\Phi_{21}(t,\omega)&=\int^{t}_{0} \Phi_{22}(t-s,\theta_s\omega)\Phi_{11}(s,\omega)ds\\
		&\leq \int^{t}_{0} e^{-(n-1)\lambda (t-s)}R_{22}(t-s, \theta_s\omega)e^{-n\lambda s}R_{11}(s, \omega)ds\\
		&=e^{-(n-1)\lambda t}\int^{t}_{0}e^{-\lambda s}{R_{22}(t-s,\theta_{s}\omega)R_{11}(s,\omega)ds}\\
		&=e^{-(n-1)\lambda t}R_{21}(t,\omega) 
	\end{align*}
	and 
	\begin{align*}
		\Phi_{31}(t,\omega)&=\int^{t}_{0} \Phi_{33}(t-s,\theta_s\omega)\Phi_{21}(s,\omega)ds\\
		&\leq \int^{t}_{0} e^{-(n-2)\lambda (t-s)}R_{33}(t-s, \theta_s\omega)e^{-(n-1)\lambda s}R_{21}(s, \omega)ds\\
		&=e^{-(n-2)\lambda t}\int^{t}_{0}e^{-\lambda s}{R_{33}(t-s,\theta_{s}\omega)R_{21}(s,\omega)ds} \\
		&=e^{-(n-2)\lambda t}R_{31}(t,\omega) 
	\end{align*}
	for all $t\geq 0$ and $\omega\in\Omega$.
	
	Furthermore, we only prove that $\sup_{s\in\mathbb{R}_+}R_{21}(s,\omega)$ and $\sup_{s\in\mathbb{R}_+}R_{31}(s,\omega)$ are tempered random variables, the rest can be analogously completed by induction. For any $\gamma>0$, by (\ref{57}) we have the following estimate
	\begin{align*}
		&~~~~\sup_{t\in\mathbb{R}}\big\{e^{-\gamma |t|}\sup_{s\in\mathbb{R_+}}R_{21}(s,\theta_t\omega)\big\}\\
		&=\sup_{s\in\mathbb{R_+},t\in\mathbb{R}}\Big\{e^{-\gamma|t|}\displaystyle \int^{s}_{0}e^{-\lambda r}{R_{22}(s-r,\theta_{r+t}\omega)R_{11}(r,\theta_t\omega)dr}\Big\}\\
		&\leq \sup_{s\in\mathbb{R_+},t\in\mathbb{R}}\Big\{\displaystyle \int^{s}_{0}e^{-\frac{1}{2}\lambda r}e^{-\frac{\lambda\wedge \gamma}{2}(|t|+r)}{R_{22}(s-r,\theta_{r+t}\omega)e^{-\frac{\gamma|t|}{2}}R_{11}(r,\theta_t\omega)dr}\Big\}\\
		&\leq \sup_{t\in\mathbb{R}}\Big\{\displaystyle \int^{\infty}_{0}e^{-\frac{1}{2}\lambda r}e^{-\frac{\lambda\wedge \gamma}{2}(|t+r|)}{\sup_{s\in\mathbb{R_+}}R_{22}(s,\theta_{r+t}\omega)e^{-\frac{\gamma|t|}{2}}\sup_{s\in\mathbb{R_+}}R_{11}(s,\theta_t\omega)dr}\Big\}\\
		&\leq \sup_{t\in\mathbb{R}}\big\{e^{-\frac{(\lambda\wedge \gamma)|t|}{2}}\sup_{s\in\mathbb{R_+}}R_{22}(s,\theta_{t}\omega)\big\}\sup_{t\in\mathbb{R}}\big\{e^{-\frac{\gamma|t|}{2}}\sup_{s\in\mathbb{R_+}}R_{11}(s,\theta_t\omega)\big\}\displaystyle \int^{\infty}_{0}e^{-\frac{1}{2}\lambda r}dr\\
		&<\infty,
	\end{align*}
	and  
	\begin{align*}
		&~~~~\sup_{t\in\mathbb{R}}\big\{e^{-\gamma |t|}\sup_{s\in\mathbb{R_+}}R_{31}(s,\theta_t\omega)\big\}\\
		&=\sup_{s\in\mathbb{R_+},t\in\mathbb{R}}\Big\{e^{-\gamma|t|}\displaystyle \int^{s}_{0}e^{-\lambda r}{R_{33}(s-r,\theta_{r+t}\omega)R_{21}(r,\theta_t\omega)dr}\Big\}\\
		&\leq \sup_{s\in\mathbb{R_+},t\in\mathbb{R}}\Big\{\displaystyle \int^{s}_{0}e^{-\frac{1}{2}\lambda r}e^{-\frac{\lambda\wedge \gamma}{2}(|t|+r)}{R_{33}(s-r,\theta_{r+t}\omega)e^{-\frac{\gamma|t|}{2}}R_{21}(r,\theta_t\omega)dr}\Big\}\\
		&\leq \sup_{t\in\mathbb{R}}\big\{e^{-\frac{(\lambda\wedge \gamma)|t|}{2}}\sup_{s\in\mathbb{R_+}}R_{33}(s,\theta_{t}\omega)\big\}\sup_{t\in\mathbb{R}}\big\{e^{-\frac{\gamma|t|}{2}}\sup_{s\in\mathbb{R_+}}R_{21}(s,\theta_t\omega)\big\}\displaystyle \int^{\infty}_{0}e^{-\frac{1}{2}\lambda r}dr\\
		&<\infty.
	\end{align*}
	Let
	\begin{equation}\label{58}
		R(t,\omega):=\max _{i,j=1, \cdots, n}R_{ij}(t,\omega).
	\end{equation}
	Then $\sup_{t\in\mathbb{R}_+}R(t,\omega)$ is a tempered random variable and
	\begin{equation*}
		\parallel \Phi(t,\omega) \parallel:= \max \{|\Phi_{ij}(t,\omega)| :i,j=1,...,n\}\leq R(t,\omega)e^{-\lambda t}, \ \ t\geq 0, \ \omega\in\Omega.
	\end{equation*}
	It follows from (\ref{22}), (\ref{55}) and (\ref{57}) that $R_{ij}(t-s,\theta_s\omega)$, $i,j=1, \cdots, n$ are ${\cal F}_s^t$-measurable, and thus $R(t-s, \theta_s\omega)$ is ${\cal F}_s^t$-measurable.
	By (\ref{55}) and the maximal inequality of geometric Brownian motion (see \cite{GP}, p.858 and \cite{PG}, p.1639), we obtain that
	\begin{equation}\label{59}
		\mathbb{E}(\sup_{t\in\mathbb{R}_+}R_{ii}(t,\omega))\leq 1+\frac{\sigma_i^2}{2i\lambda}, \ \ i=1, \cdots, n.
	\end{equation}
	\quad Combining (\ref{57}), (\ref{59}), the $\mathbb{P}$-measure preserving property of $\theta$, Fubini theorem and the fact that an $n$-dimensional Brownian motion has $n$ independent components, we obtain that
	\begin{align*}
		&~~~~\sup_{t\in\mathbb{R}_+}\mathbb{E}(R_{ij}(t,\omega))\\
		&=\sup_{t\in\mathbb{R}_+}\mathbb{E}\big(\int^{t}_{0}e^{-\lambda s}{R_{ii}(t-s,\theta_{s}\omega)R_{i-1, j}(s,\omega)ds}\big)\\
		&=\sup_{t\in\mathbb{R}_+}\int^{t}_{0}e^{-\lambda s}\mathbb{E}\big(R_{ii}(t-s,\theta_{s}\omega)R_{i-1, j}(s,\omega)\big)ds\\
		&=\sup_{t\in\mathbb{R}_+}\int^{t}_{0}e^{-\lambda s}\mathbb{E}(R_{ii}(t-s,\theta_{s}\omega))\mathbb{E}(R_{i-1, j}(s,\omega))ds\\
		&\leq\sup_{t\in\mathbb{R}_+}\mathbb{E}(R_{i-1, j}(t,\omega))\sup_{t\in\mathbb{R}_+}\int^{t}_{0}e^{-\lambda s}\mathbb{E}(R_{ii}(t-s,\omega))ds\nonumber\\
		&\leq \sup_{t\in\mathbb{R}_+}\mathbb{E}(R_{ii}(t,\omega))\sup_{t\in\mathbb{R}_+}\mathbb{E}(R_{i-1, j}(t,\omega))\int^{\infty}_{0}e^{-\lambda s}ds\nonumber\\
		&\leq \frac{1}{\lambda}\sup_{t\in\mathbb{R}_+}\mathbb{E}(R_{ii}(t,\omega))\sup_{t\in\mathbb{R}_+}\mathbb{E}(R_{i-1, j}(t,\omega)),\ \ 1\leq j\leq i-1.
	\end{align*}
	By induction, we conclude that
	\begin{equation*}
		\sup_{t\in\mathbb{R}_+}\mathbb{E}(R_{ij}(t,\omega))\leq\frac{1}{\lambda^{i-j}}\prod_{k=j}^{i}\sup_{t\in\mathbb{R}_+}\mathbb{E}(R_{kk}(t,\omega)), \ \ 1\leq j\leq i.
	\end{equation*}
	So 
	\begin{align*}
		\sup_{t\in\mathbb{R}_+}\mathbb{E}(R(t,\omega))&\leq\sup_{t\in\mathbb{R}_+}\sum_{i=1}^{n}\sum_{j=1}^{i}\mathbb{E}(R_{ij}(t,\omega))\nonumber\\
		&\leq\sum_{i=1}^{n}\sum_{j=1}^{i}\sup_{t\in\mathbb{R}_+}\mathbb{E}(R_{ij}(t,\omega))\nonumber\\
		&\leq\sum_{i=1}^{n}\sum_{j=1}^{i}\frac{1}{\lambda^{i-j}}\prod_{k=j}^{i}\sup_{t\in\mathbb{R}_+}\mathbb{E}(R_{kk}(t,\omega))\nonumber\\
		&\leq\sum_{i=1}^{n}\sum_{j=1}^{i}\frac{1}{\lambda^{i-j}}\prod_{k=j}^{i}\mathbb{E}(\sup_{t\in\mathbb{R}_+}R_{kk}(t,\omega))\nonumber\\
		&\leq \sum_{i=1}^{n}\sum_{j=1}^{i}\frac{1}{\lambda^{i-j}}\prod_{k=j}^{i}(1+\frac{\sigma_k^2}{2k\lambda}).
	\end{align*}
	
	Let $h(t, x)=(f(t, x_n), 0, \cdots, 0)^T$, $x\in\mathbb{R}_+^n$, $N_1=\sup_{t\in\mathbb{R}, x_n\in\mathbb{R}_+}|f(t, x_n)|$, $N_i=0$ for all $2\leq i\leq n$ and $L=\sup_{t\in\mathbb{R}, x_n\in\mathbb{R}_+}\left|\frac{df(t,x_n)}{d x_n}\right|$. Then from Theorem \ref{4.3} we conclude the following.
	\begin{pro}\label{5.1}
		Let $\alpha_i>0$ for $i=1,\cdots,n$ and $f\in C_b^1(\mathbb{R}\times\mathbb{R}_+, \mathbb{R}_+)$. Assume that $f(t+T, x)=f(t, x)$, $t\in\mathbb{R}$, $x\in\mathbb{R_+}$ and $f(t, \cdot)$, $t\in\mathbb{R}$ is increasing or decreasing in $\mathbb{R}_+$. If 
		\begin{equation*}
			\frac{Ln^2 \sup_{t\in\mathbb{R}_+}\mathbb{E}(R(t,\omega))}{\lambda}\leq\frac{Ln^2}{\lambda}\sum_{i=1}^{n}\sum_{j=1}^{i}\frac{1}{\lambda^{i-j}}\prod_{k=j}^{i}(1+\frac{\sigma_k^2}{2k\lambda})<1
		\end{equation*}
		holds, where $R(t, \omega)$ is defined in (\ref{58}) and $L=\sup_{t\in\mathbb{R}, x_n\in\mathbb{R}_+}\left|\frac{df(t,x_n)}{d x_n}\right|$. Then equation (\ref{51}) has a unique random periodic solution of periodic T in $\mathbb{R}_+^n$.
	\end{pro}
	
	\begin{exa}(Time-periodic stochastic Goodwin system). Consider $n$-dimensional stochastic differential equation
		\begin{eqnarray}\label{510}
			\left\{
			\begin{array}{l}
				dx_1=(-\alpha_1 x_1+\frac{V}{K+\sin t+x_n^m})dt+\sigma_1x_1dW_t^1,\\
				dx_i=(x_{i-1}-\alpha_i x_i)dt+\sigma_ix_idW_t^i, \ 2\leq i\leq n.
			\end{array}
			\right.
		\end{eqnarray}
		where $m>1$, $K>2$, $V>0$ and $\alpha_i>0$ for $i=1,\cdots,n$. It is clear that (\ref{510}) is a non-monotone system, which can be regarded as the stochastic Goodwin model with time-periodic coefficient; see \cite{G, HTW}. By the direct calculation, it is obvious that
		\begin{equation*}
			L=\sup_{t\in\mathbb{R}, x_n\in\mathbb{R}_+}\left|\frac{mVx_n^{m-1}}{(K+\sin t+x_n^m)^2}\right|\leq \sup_{t\in\mathbb{R}, x_n\in\mathbb{R}_+}\left|\frac{mV(1+x_n^m)}{(K+\sin t+x_n^m)^2}\right|\leq \frac{mV}{K-1}.
		\end{equation*}
		Applying Proposition \ref{5.1}, we get that if
		\begin{equation}\label{511}
			\frac{mn^2V}{\lambda(K-1)}\sum_{i=1}^{n}\sum_{j=1}^{i}\frac{1}{\lambda^{i-j}}\prod_{k=j}^{i}(1+\frac{\sigma_k^2}{2k\lambda})<1
		\end{equation}
		is satisfied, then (\ref{510}) has a unique globally stable random periodic solution of periodic $2\pi$ in $\mathbb{R}_+^n$. Here, (\ref{511}) holds for $V$ sufficiently small or $K$ sufficiently large. 
	\end{exa}
	
	\begin{exa} (Time-periodic stochastic Othmer-Tyson system). Consider the following n-dimensional nonautonomous stochastic Othmer-Tyson positive feedback system:
		\begin{eqnarray}\label{512}
			\left\{
			\begin{array}{l}
				dx_1=(-\alpha_1 x_1+\frac{k_0(1+x_n^m)}{K+\sin t+x_n^m})dt+\sigma_1x_1dW_t^1,\\
				dx_i=(x_{i-1}-\alpha_i x_i)dt+\sigma_ix_idW_t^i, \ 2\leq i\leq n.
			\end{array}
			\right.
		\end{eqnarray}
		where $m>1$, $K>2$, $k_0>0$ and $\alpha_i>0$ for $i=1,\cdots,n$. This is a stochastic cooperative system; see \cite{OH, TO}. By the direct calculation, it is obvious that
		\begin{equation*}
			L=\sup_{t\in\mathbb{R}, x_n\in\mathbb{R}_+}\left|\frac{mk_0(K+\sin t-1)x_n^{m-1}}{(K+\sin t+x_n^m)^2}\right|\leq \sup_{t\in\mathbb{R}, x_n\in\mathbb{R}_+}\left|\frac{mk_0(K+\sin t-1)(1+x_n^m)}{(K+\sin t+x_n^m)^2}\right|\leq \frac{mk_0K}{K-1}.
		\end{equation*}
		Applying Proposition \ref{5.1}, we get that if
		\begin{equation}\label{513}
			\frac{mk_0n^2K}{\lambda(K-1)}\sum_{i=1}^{n}\sum_{j=1}^{i}\frac{1}{\lambda^{i-j}}\prod_{k=j}^{i}(1+\frac{\sigma_k^2}{2k\lambda})<1
		\end{equation}
		is satisfied, then (\ref{512}) has a unique globally stable random periodic solution of periodic $2\pi$ in $\mathbb{R}_+^n$. Here, (\ref{513}) holds for $k_0$ sufficiently small.
	\end{exa}

	\begin{exa}
		Consider an $n$-dimensional stochastic competitive system with time-periodic coefficient
		\begin{equation}\label{514}
			dx_i=[-\alpha_ix_i+h_i(t,x)]dt+\sigma_ix_idW_t^i, 
		\end{equation}
		where $\alpha_i>0$ for all $i=1, \cdots, n$ and  
		\begin{equation*}
			h_i(t,x):=|\sin t|+\frac{1}{K_i+x_1^m+\cdots+x_n^m}, \ \ x\in\mathbb{R}_+^n, \ i=1, \cdots, n,
		\end{equation*}
		where $m>1$ and $K_i>1$ for all $i=1, \cdots, n$. Then $h(t, \cdot)$ is a decreasing function from $\mathbb{R}_+^n$ to $\mathbb{R}_+^n\setminus 0$. The fundamental matrix of the following corresponding linear homogeneous stochastic It$\hat{o}$ type differential equations
		\begin{equation*}
			dx_i=-\alpha_ix_idt+\sigma_ix_idW_t^i, 
		\end{equation*}
		is
		$$
		\Phi(t, \omega)=\begin{pmatrix}
			\Phi_{11}(t,\omega) & 0 & \cdots & 0 \\
			0 & \Phi_{22}(t,\omega) & \cdots & 0 \\
			\vdots & \vdots & \ddots & \vdots \\
			0 & 0 & \cdots & \Phi_{nn}(t,\omega) \\ 
		\end{pmatrix}
		$$
		for all $t\geq 0$ and $\omega\in\Omega$, where
		\begin{equation*}
			\Phi_{ii}(t, \omega)= e^{-(\alpha_i+\frac{1}{2}\sigma_i^2)t+\sigma_i W_t^i(\omega)}.
		\end{equation*}
		It is clear that
		\begin{equation*}
			\parallel \Phi(t,\omega) \parallel:= \max \{|\Phi_{ij}(t,\omega)| :i,j=1,...,n\}\leq R(t,\omega)e^{-\lambda t}, \ \ t\geq 0, \ \omega\in\Omega,
		\end{equation*}
		where $\lambda=\frac{1}{2}\min\{\alpha_1, \cdots, \alpha_n\}$ and 
		\begin{equation*}
			R(t, \omega)=\bigvee_{i=1}^ne^{-(\lambda+\frac{1}{2}\sigma_i^2)t+\sigma_i W_t^i(\omega)}.
		\end{equation*}
		Similar to the analysis for system (\ref{51}), it is easy to check that $R(t-s,\theta_s\omega)$ is ${\cal F}_s^t$-measurable and for any $\gamma>0$,
		\begin{equation*}
			\sup_{t\in\mathbb{R}}\big\{e^{-\gamma |t|}\sup_{s\in\mathbb{R_+}}R(s,\theta_t\omega)\big\}<\infty.
		\end{equation*}
		By the maximal inequality of geometric Brownian motion, we can obtain that
		\begin{align*}
			\sup_{t\in\mathbb{R}_+}\mathbb{E}(R(t,\omega))&=\sup_{t\in\mathbb{R}_+}\mathbb{E}\Big(\bigvee_{i=1}^ne^{-(\lambda+\frac{1}{2}\sigma_i^2)t+\sigma_i W_t^i(\omega)}\Big)\\
			&\leq\sum_{i=1}^{n}\sup_{t\in\mathbb{R}_+}\mathbb{E}e^{-(\lambda+\frac{1}{2}\sigma_i^2)t+\sigma_i W_t^i(\omega)}\nonumber\\
			&\leq\sum_{i=1}^{n}\mathbb{E}(\sup_{t\in\mathbb{R}_+}e^{-(\lambda+\frac{1}{2}\sigma_i^2)t+\sigma_i W_t^i(\omega)})\nonumber\\
			&\leq \sum_{i=1}^{n}(1+\frac{\sigma_i^2}{2\lambda}).
		\end{align*}
		By the direct calculation, it is obvious that
		\begin{align*}
			L&=\max \left\{\ \sup_{t\in\mathbb{R}, x\in\mathbb{R}_+^d} \left|\frac{\partial h_i(t,x)}{\partial x_j}\right|, i,j=1,...,n\right\}\\
			&=\max \left\{\ \sup_{x\in\mathbb{R}_+^d} \left|\frac{mx_j^{m-1}}{(K_i+x_1^m+\cdots+x_n^m)^2}\right|, i,j=1,...,n\right\}\\
			&\leq \frac{m}{K_i}\leq\frac{m}{K},
		\end{align*}
		where $K=\min\{K_i, i=1, \cdots, n\}$. Therefore, if
		\begin{equation}\label{515}
			\frac{mn^2}{\lambda K}\sum_{i=1}^{n}(1+\frac{\sigma_i^2}{2\lambda})<1
		\end{equation}
		is satisfied, then (\ref{514}) has a unique globally stable random periodic solution of periodic $\pi$ in $\mathbb{R}_+^n$. Here, (\ref{515}) holds when $\lambda$ or $K$ is large enough.
	\end{exa} 
		As a further specific example, we give the following example.
	\begin{exa}
	 Consider the following 3-dimensional stochastic Othmer-Tyson positive feedback system with time-periodic coefficient:
		\begin{eqnarray}\label{516}
			\left\{
			\begin{array}{l}
				dx_1=(-8 x_1+\frac{1}{12}\cdot \frac{1+x_3^3}{3+\sin t+x_3^3})dt+\frac{1}{2}x_1dW_t^1,\\
				dx_2=(x_1-9x_2)dt+\frac{1}{4}x_2dW_t^2, \\
				dx_3=(x_2-10x_3)dt+\frac{1}{3}x_3dW_t^3.
			\end{array}
			\right.
		\end{eqnarray}
		By the direct calculation, it is obvious that
		\begin{align*}
			L&=\sup_{t\in\mathbb{R}, x_3\in\mathbb{R}_+}\left|\frac{(2+\sin t)x_3^2}{4(3+\sin t+x_3^3)^2}\right|\\
			&=\sup_{t\in\mathbb{R}} \frac{(2+\sin t)x_3^2}{4(3+\sin t+x_3^3)^2}\Big|_{x_3^3= \frac{3+sint}{2}}\\
			&=\sup_{t\in\mathbb{R}}
			{\frac{2+\sin t}{36(\frac{3+\sin t}{2})^\frac{4}{3}}}\\
			&=\frac{2+\sin t}{36(\frac{3+\sin t}{2})^\frac{4}{3}}\Big|_{\sin t=1}\\
			&=\frac{1}{24\cdot 2^{\frac{1}{3}}}.
		\end{align*}
		
		The corresponding linear homogeneous stochastic It$\hat{o}$ type differential equations is
		\begin{eqnarray}\label{517}
			\left\{
			\begin{array}{l}
				dx_1=-8 x_1dt+\frac{1}{2}x_1dW_t^1,\\
				dx_2=(x_1-9x_2)dt+\frac{1}{4}x_2dW_t^2, \\
				dx_3=(x_2-10x_3)dt+\frac{1}{3}x_3dW_t^3.
			\end{array}
			\right.
		\end{eqnarray}
		By the variation-of-constants formula, we can easily calculate the fundamental matrix $\Phi(t,\omega)$ of (\ref{517}) as follows:
		$$
		\Phi(t, \omega)=\begin{pmatrix}
			\Phi_{11}(t,\omega) & 0 & 0 \\
			\Phi_{21}(t,\omega) & \Phi_{22}(t,\omega) & 0 \\
			\Phi_{31}(t,\omega) & \Phi_{32}(t,\omega) & \Phi_{33}(t,\omega) \\ 
		\end{pmatrix}
		$$
		for all $t\geq 0$ and $\omega\in\Omega$, where
		\begin{equation*}
			\Phi_{11}(t,\omega)= e^{-(8+\frac{1}{8})t+\frac{1}{2} W_t^1(\omega)},
		\end{equation*}
		\begin{equation*}
			\Phi_{22}(t,\omega)= e^{-(9+\frac{1}{32})t+\frac{1}{4} W_t^2(\omega)},
		\end{equation*}
		\begin{equation*}
			\Phi_{33}(t,\omega)= e^{-(10+\frac{1}{18})t+\frac{1}{3} W_t^3(\omega)}.
		\end{equation*}
		and 
		\begin{equation*}
			\Phi_{ij}(t, \omega)=\left\{
			\begin{array}{rcl}
				\displaystyle \int^{t}_{0}{\Phi_{ii}(t-s,\theta_{s}\omega)\Phi_{i-1, j}(s,\omega))ds}, & & {1\leq j \leq i-1}\\ 
				0, & & {i+1\leq j\leq 3.}
			\end{array} \right.
		\end{equation*}
		Hence, it is easy to check that
		\begin{equation*}
			\Phi_{ii}(t, \omega)= R_{ii}(t,\omega)e^{-(4-i)t}, \ \ i=1, 2, 3,
		\end{equation*}
		for all $t\geq 0$ and $\omega\in\Omega$, where
		\begin{equation*}
			R_{11}(t, \omega)=e^{-(5+\frac{1}{8})t+\frac{1}{2} W_t^1(\omega)}, 
		\end{equation*}
		\begin{equation*}
			R_{22}(t, \omega)=e^{-(7+\frac{1}{32})t+\frac{1}{4} W_t^2(\omega)}, 
		\end{equation*}		
		\begin{equation*}
			R_{33}(t, \omega)=e^{-(9+\frac{1}{18})t+\frac{1}{3} W_t^3(\omega)}. 
		\end{equation*}		
		and thus
		\begin{align*}
			\Phi_{21}(t, \omega)&=\int^{t}_{0}{\Phi_{22}(t-s,\theta_{s}\omega)\Phi_{11}(s,\omega))ds}\\
			&=\int^{t}_{0}{e^{-2(t-s)}R_{22}(t-s,\theta_{s}\omega)e^{-3s}R_{11}(s,\omega))ds}\\
			&=e^{-2t}R_{21}(t,\omega),
		\end{align*}
		where $ R_{21}(t,\omega)=\int^{t}_{0}{e^{-s}R_{22}(t-s,\theta_{s}\omega)R_{11}(s,\omega))ds}$.
		\begin{align*}
			\Phi_{31}(t, \omega)&=\int^{t}_{0}{\Phi_{33}(t-s,\theta_{s}\omega)\Phi_{21}(s,\omega))ds}\\
			&=\int^{t}_{0}{e^{-(t-s)}R_{33}(t-s,\theta_{s}\omega)e^{-2s}R_{21}(s,\omega))ds}\\
			&=e^{-t}R_{31}(t,\omega),
		\end{align*}
		where $R_{31}(t,\omega)=\int^{t}_{0}{e^{-s}R_{33}(t-s,\theta_{s}\omega)R_{21}(s,\omega))ds}$,	
		\begin{align*}
			\Phi_{32}(t, \omega)&=\int^{t}_{0}{\Phi_{33}(t-s,\theta_{s}\omega)\Phi_{22}(s,\omega))ds}\\
			&=\int^{t}_{0}{e^{-(t-s)}R_{33}(t-s,\theta_{s}\omega)e^{-2s}R_{22}(s,\omega))ds}\\
			&=e^{-t}R_{32}(t,\omega),
		\end{align*}	
		where $R_{32}(t,\omega)=\int^{t}_{0}{e^{-s}R_{33}(t-s,\theta_{s}\omega)R_{22}(s,\omega))ds}$.	
		In order to verify conditions in Theorem \ref{4.3}, we choose $\lambda=1$ and 
		$$R(t,\omega)=\max _{i,j=1, 2, 3}R_{ij}(t,\omega).$$
		Similar to the analysis for system (\ref{51}), we can prove that $\sup_{t\in\mathbb{R}_+}R(t,\omega)$ is a tempered random variable and
		\begin{equation*}
			\parallel \Phi(t,\omega) \parallel:= \max \{|\Phi_{ij}(t,\omega)| :i,j=1, 2, 3\}\leq R(t,\omega)e^{-t}, \ \ t\geq 0, \ \omega\in\Omega.
		\end{equation*}
		It is clear that $R_{ij}(t-s,\theta_s\omega)$, $i,j=1, \cdots, n$ is ${\cal F}_s^t$-measurable, and thus $R(t-s, \theta_s\omega)$ is ${\cal F}_s^t$-measurable. By the maximal inequality of geometric Brownian motion, i.e.,  $$\mathbb{E}\left(\sup_{t\in\mathbb{R}_+}e^{(\mu-\frac{1}{2}\sigma^2)t+\sigma W_t(\omega)}\right)=1-\frac{\sigma^2}{2\mu}$$
		 for $\mu<0$ and $\sigma>0$, we obtain that
		\begin{align*}
			\sup_{t\in\mathbb{R}_+}\mathbb{E}(R(t,\omega))&=\sup_{t\in\mathbb{R}_+}\mathbb{E}(\max _{i,j=1, 2, 3}R_{ij}(t,\omega))\nonumber\\
			&\leq\mathbb{E}(\sup_{t\in\mathbb{R}_+}R_{31}(t,\omega))+\mathbb{E}(\sup_{t\in\mathbb{R}_+}R_{32}(t,\omega))+\mathbb{E}(\sup_{t\in\mathbb{R}_+}R_{33}(t,\omega))\nonumber\\
			&\leq\frac{41}{40}\times\frac{225}{224}\times\frac{163}{162}+\frac{225}{224}\times\frac{163}{162}+\frac{163}{162} \nonumber\\
			&<3.0528.
		\end{align*}
		Therefore, we have that 
		$$ \frac{Ld^2 \sup_{t\in\mathbb{R}_+}\mathbb{E}(R(t,\omega))}{\lambda}\leq \frac{1}{24\times 2^{\frac{1}{3}}}\times 9 \times 3.0528 <1.$$	
		Then by Theorem \ref{4.3}, we conclude that (\ref{516}) has a unique globally stable random periodic solution of period $2\pi$ in $\mathbb{R}_+^3$.
	\end{exa}
	
	{ \noindent {\bf\large Acknowledgements}\ \
		This work was supported in part by National Key R\&D Program of China (No. 2020YFA0712700), National Natural Science Foundation of China (Nos. 11931004, 12031020, 12090014), CAS Key Project of Frontier Sciences (No. QYZDJ-SSW-JSC003), the Key Laboratory of Random Complex Structures and Data Sciences, CAS(No. 2008DP173182).

	\end{document}